\documentclass[a4paper,11pt]{amsart}
\usepackage{amssymb}
\newtheorem{lemma}{Lemma}
\newtheorem{rem}{Remark}
\newtheorem{defi}{Definition}
\newtheorem{prop}{Proposition}
\newtheorem*{thm}{Theorem}
\newtheorem{theo}{Theorem}
\newtheorem{cor}{Corollary}

\newenvironment{bulletlist}{\begin{list}{\labelitemi}%
{\setlength{\leftmargin}{\parindent}\def
\makelabel ##1{\hss \llap {\upshape ##1}}}}{\end{list}}
\newcommand{\R}{{\mathbb R}}
\newcommand{\C}{{\mathbb C}}
\newcommand{\E}{{\mathbb E}}
\newcommand{\Ric}{\mathit{Ric}}
\newcommand{\Scal}{\mathit{Scal}}
\newcommand{\Id}{\mathit{Id}}
\newcommand{\g}[1]{\langle #1 \rangle}
\newcommand{\ka}{K{\"a}hler }
\def\biax/{ortho-toric}
\def\monax/{mono-axial}
\def\biham/{hamiltonian}
\def\monoham/{hamiltonian}
\begin{document}
\title[Weakly selfdual K\smash{\"a}hler surfaces]{The geometry of weakly
selfdual\\ K{\"a}hler surfaces}
\author{V. Apostolov, D. M. J. Calderbank and P. Gauduchon}
\thanks{The first author was supported in part by a PAFARC-UQAM grant,
and by an NSERC grant, the second author by the Leverhulme Trust and the
William Gordon Seggie Brown Trust. All three authors are members of EDGE,
Research Training Network HPRN-CT-2000-00101, supported by the European
Human Potential Programme.}
\address{Vestislav Apostolov \\ D{\'e}partement de Math{\'e}matiques\\
UQAM\\ C.P. 8888 \\ Succ. Centre-ville \\ Montr{\'e}al (Qu{\'e}bec) \\
H3C 3P8 \\ Canada}
\email{apostolo@math.uqam.ca}
\address{David M. J. Calderbank \\ Department of Mathematics and
  Statistics\\ University of Edinburgh\\ King's Building\\ Mayfield's
  Road\\ Edinburgh EH9 3JZ\\ Scotland}
\email{davidmjc@maths.ed.ac.uk} 
\address{Paul Gauduchon \\ Centre de Math\'ematiques\\
{\'E}cole Polytechnique \\ UMR 7640 du CNRS
\\ 91128 Palaiseau \\ France}
\email{pg@math.polytechnique.fr}
\date{April 2001}
\begin{abstract} We study K{\"a}hler surfaces with harmonic
anti-selfdual Weyl tensor. We provide an explicit
local description, which we use to obtain the complete classification
in the compact case. We give new examples of extremal K{\"a}hler
metrics, including K{\"a}hler--Einstein metrics and conformally Einstein
K{\"a}hler metrics. We also extend some of our results to almost K{\"a}hler
$4$-manifolds, providing new examples of Ricci-flat almost K{\"a}hler
metrics which are not K{\"a}hler.
\end{abstract}
\maketitle

\section*{Introduction}

Selfdual K{\"a}hler surfaces have been considered in several recent works,
in particular in a paper by R. Bryant \cite{bryant}, where selfdual
K{\"a}hler surfaces appear as the four-dimensional case of a comprehensive
study of Bochner-flat K{\"a}hler manifolds in all dimensions, and in a paper
by two of the authors \cite{AG2}, where a generic equivalence has been
established between selfdual K{\"a}hler surfaces and selfdual Hermitian
Einstein metrics and where an explicit local description of the latter is
provided.

Whereas selfdual surfaces are easily proved to be extremal, i.e., admitting
a hamiltonian Killing vector field whose momentum map is the scalar
curvature, it was an a priori unexpected fact, independently discovered in
the above works, that they actually admit a second hamiltonian Killing
vector field; moreover, a crucial observation of R. Bryant \cite{bryant}, is
that the momentum map of the latter is the pfaffian of the normalized Ricci
form. Since these Killing vector fields commute, this also provides a link
with the work of H.~Pedersen and the second author~\cite{CaPe:sdt2}, where
an explicit local classification of selfdual Einstein metrics with two
commuting Killing vector fields is obtained, without the hypothesis that
they are Hermitian.

In this paper, we show that we can relax the assumption of selfduality, and
establish the same bi-hamiltonian structure for {\it weakly selfdual}
K{\"a}hler surfaces, i.e., K{\"a}hler surfaces whose anti-selfdual Weyl
tensor $W ^-$ is harmonic; this fact has its origins in the basic {\it
Matsumoto--Tanno identity} for such surfaces, recently re-discovered by
W. Jelonek \cite{jelonek}, and leads to a surprisingly simple explicit
expression, generalizing an expression found by Bryant in the selfdual case.
On the other hand, we also observe that in Calabi's family of extremal
K{\"a}hler metrics on the first Hirzebruch surface $F_1$, there is a unique
(and completely explicit) weakly selfdual metric up to homothety. This is
in contrast to the selfdual case, where the compact (smooth) examples are
all locally symmetric~\cite{BYC}.

Our results concerning weakly selfdual K{\"a}hler surfaces may be summarized
as follows (definitions and more precise statements are given in the body of
the paper).

\begin{thm} Let $(M,g,J,\omega)$ be a weakly selfdual K{\"a}hler surface.
Then $(g,J)$ is a \emph{bi-extremal K{\"a}hler metric} in the sense that
the scalar curvature and the pfaffian of the normalized Ricci form of
$(g,J)$ are Poisson-commuting momentum maps for hamiltonian Killing vector
fields $K_1$ and $K_2$ respectively. Furthermore, on each connected
component of $M$ one of the following holds.

\smallskip\noindent{\rm (i)}
$K_1$ and $K_2$ are linearly independent on a dense open set. Then
$(g,J,\omega)$ has an explicit local
form~\eqref{eq:newgbis}--\eqref{eq:Gter}, depending on an arbitrary
polynomial of degree $4$ and an arbitrary constant which is zero if and only
if $g$ is selfdual, cf.~Theorem {\rm 2}.

\smallskip\noindent{\rm (ii)}
$K_1$ is non-vanishing on a dense open set, but $K_1\wedge K_2$
is identically zero. Then $(g,J,\omega)$ is locally of cohomogeneity one and
is given explicitly by the Calabi construction, cf.~Theorem {\rm 4}.

\smallskip\noindent{\rm (iii)}
$K_1$ and $K_2$ vanish identically. Then $g$ has parallel Ricci curvature
\textup(hence is either K{\"a}hler--Einstein or locally a K{\"a}hler product
of two Riemann surfaces of constant curvatures\textup).

\smallskip If $(M,g,J,\omega)$ is compact and connected, then it necessarily
belongs to case {\rm (ii)} or {\rm (iii)} above, and in case {\rm (ii)}
$(M,g,J,\omega)$ is isomorphic to the weakly selfdual Calabi extremal metric
on $F_1$ \textup(see Theorem {\rm 5}\textup).
\end{thm}

The path to proving this result touches upon various important themes in
K{\"a}hler geometry, and along the way we introduce further ideas, results
and examples. A key aspect of our approach is to study weakly selfdual
K{\"a}hler surfaces within a more general setting. The Matsumoto--Tanno
identity for weakly selfdual K{\"a}hler surfaces is equivalent to the fact
that the primitive part $\rho_0$ of the Ricci form $\rho$ of $(g,J)$
satisfies an overdetermined linear differential equation. On the open set
where $\rho_0$ is nonvanishing, the equation means that $\rho_0$ defines a
conformally K{\"a}hler Hermitian structure $I$ inducing the opposite
orientation to $J$.  Many of the properties of weakly selfdual K{\"a}hler
surfaces are simple consequences of the fact that the $\rho$ is a closed
$J$-invariant $2$-form, whose primitive part satisfies this equation. In
particular, we prove in Proposition~\ref{prop:biham} that two of the
algebraic invariants (essentially the trace and the pfaffian) of \emph{any}
such $2$-form $\varphi$ are Poisson-commuting momentum maps for hamiltonian
Killing vector fields.  Therefore, throughout the work, we develop the
theory of K{\"a}hler surfaces with such `\biham/' $2$-forms $\varphi$, which
include other interesting examples in addition to weakly selfdual K{\"a}hler
surfaces.

First of all we study the generic case that the hamiltonian Killing vector
fields are linearly independent. This means that the K{\"a}hler structure
$(g,J,\omega)$ is toric, and in Theorem 1 we characterize the class of toric
K{\"a}hler structures arising in this way from \biham/ $2$-forms. Whereas
toric K{\"a}hler surfaces in general depend essentially on an arbitrary
function of two variables~\cite{Abreu,G:kstv}, our toric surfaces, which
we call `\biax/', have an explicit form---given in
Proposition~\ref{prop:biax}---depending only on two arbitrary functions of one
variable. This has the practical advantage that curvature conditions lead to
ordinary differential equations for these functions. In particular, we are
able to obtain explicitly all of the extremal toric K{\"a}hler structures in
our class, including some new examples of K{\"a}hler metrics which are
conformally Einstein, but neither selfdual nor anti-selfdual, and also some
explicit K{\"a}hler--Einstein metrics. The weakly selfdual metrics in this
family are classified in Theorem 2.

The case that the hamiltonian Killing vector fields associated to a
hamiltonian $2$-form $\varphi$ are linearly dependent, but not both zero, is
closely related to the Calabi construction of K{\"a}hler metrics on line
bundles over a Riemann surface~\cite{calabi1}. We provide, in Theorem 3, a
geometric local characterization of these K{\"a}hler metrics: they are the
K{\"a}hler metrics $(g,J)$, with a Killing vector field $K$ such that the
almost Hermitian pair $(g,I)$, where $I$ is equal to $J$ on span of
$\{K,JK\}$ but $-J$ on the orthogonal distribution, is conformally
K{\"a}hler. Over a fixed Riemann surface $\Sigma$, the general form of these
K{\"a}hler metrics `of Calabi type' again depends essentially only on
functions of one variable from which it is easy to recover the Calabi
extremal metrics. We present these in Proposition~\ref{prop:monoextremal}:
the Riemann surface $\Sigma$ has constant curvature, and the metrics have
local cohomogeneity one under $U(2)$, $U(1,1)$ or a central extension of the
Heisenberg group Nil. The weakly selfdual Calabi extremal metrics are
classified in Theorem 4: there is a four parameter family, one of which is
globally defined on the first Hirzebruch surface. The existence of such
a metric has been independently observed by Jelonek~\cite{jelonek2}.

The proof of the above theorem is completed by classifying the compact
weakly selfdual K{\"a}hler surfaces. A partial classification for real
analytic K{\"a}hler surfaces has been recently obtained by
Jelonek~\cite{jelonek}, but we improve it in two respects: first, as
speculated by Jelonek, we are able to remove the assumption of
real-analyticity; second we prove that the only (non-product
non-K{\"a}hler--Einstein) weakly selfdual K{\"a}hler metric on a ruled
surface, is the Calabi extremal example on the first Hirzebruch surfaces
$F_1$.

The paper is organized as follows. In the first section, we establish some
basic facts concerning weakly selfdual K{\"a}hler surfaces; in particular,
using general properties of \biham/ $2$-forms, we show that weakly selfdual
K{\"a}hler surfaces are bi-extremal and that their anti-selfdual Weyl tensor
is degenerate (some facts proved in this section also appear in Jelonek's
paper \cite{jelonek}). We also present the rough classification, Proposition
6, that allows us to deduce the above theorem from Theorems 2, 4 and 5.  The
generic, toric case, and Theorem 2, are described in the second section,
whereas the third section treats the Calabi examples and Theorem 4.  The
classification of compact weakly selfdual \ka surfaces is given in section
4.

In the final section we show that some of our results generalize to the
class of almost \ka 4-manifolds $(M,g,J,\omega)$ whose Ricci tensor is
$J$-invariant.  We first observe that the Calabi construction gives rise to
new (local) examples of selfdual, Ricci-flat almost K{\"a}hler 4-manifolds
(see Example 2), which provide further local counterexamples
(cf. \cite{Arm,NuP}) to the still open Goldberg conjecture which states that
a {\it compact} Einstein almost \ka manifold must be K{\"a}hler--Einstein.
Next we consider compact almost \ka 4-manifolds with $J$-invariant Ricci
tensor which are weakly selfdual, i.e., have harmonic anti-selfdual Weyl
tensor. We show in Theorem \ref{theo6} that weak-selfduality has strong
consequences for the integrability of the corresponding almost complex
structure, providing another interesting link with the Goldberg
conjecture. As an application of this global result, we prove that a compact
almost \ka 4-manifold has constant sectional curvatures on the Lagrangian
2-planes if and only if it is a selfdual \ka surface (see Corollary
\ref{theo7}).

\section{Weakly selfdual K{\"a}hler surfaces} \label{s:wsdk}

\subsection{The Matsumoto--Tanno identity}

A {\it K{\"a}hler surface} $(M, g, J)$ is an oriented Riemannian
four-dimensional manifold equipped with a selfdual complex structure $J$,
such that $\nabla J = 0$, where $\nabla$ denotes the Levi-Civita connection
of $g$. The K{\"a}hler form is the $J$-invariant selfdual $2$-form $\omega
(\cdot, \cdot) = \g{J\cdot, \cdot}$; $\omega$ is closed and $(M, \omega)$
is a symplectic manifold.

The vector bundle $\Lambda ^+ M$ of selfdual $2$-forms is the orthogonal
direct sum of the trivial bundle generated by $\omega$ and of the bundle of
$J$-anti-invariant $2$-forms, whereas the bundle $\Lambda ^- M$ coincides
with the bundle of {\it primitive}---or {\it trace-free}---$J$-invariant
$2$-forms.

We denote by $R$, $\Ric$, $\Ric _0$, $\Scal$, $W = W ^+ + W ^-$, the
curvature, the Ricci tensor, the trace-free part of $\Ric$, the scalar
curvature (i.e., the trace of $\Ric$), and the Weyl tensor, expressed as
the sum of its $\pm$-selfdual components $W ^{\pm}$.

The {\it Ricci form}, $\rho$, of a K{\"a}hler surface is the $J$-invariant
$2$-form defined by $\rho (\cdot, \cdot) = \Ric (J \cdot, \cdot)$; $\rho$ is
closed and, up to a factor $2 \pi$, is a representative of the first Chern
class of $(M, J)$ in de Rham cohomology; the trace-free part of $\rho$ is
denoted by $\rho _0$.

\begin{defi} \label{def1} A K{\"a}hler surface $(M, g, J)$ is {\rm weakly
selfdual} if its anti-selfdual Weyl tensor $W ^-$ is harmonic, i.e.,
satisfies:
\begin{equation*} \delta ^g W ^- = 0. \end{equation*}
where the codifferential $\delta ^g$ acts on $W^-$ as on a $2$-form with
values in $\Lambda ^- M$.
\end{defi}

Because of the Bianchi identity, the weak selfduality condition can also
be defined in terms of the Ricci tensor; from this point of view, it can
also be considered as a {\it weak Einstein condition}, in the sense that
every Einstein metric is weakly selfdual. The link is provided by the {\it
Cotton--York} tensor $C_{X,Y}(Z)$ of the Riemannian metric $g$.  Recall that
the Cotton--York tensor is defined by $C _{X, Y}(Z) = - (\nabla _X h )(Y,
Z) + (\nabla_Y h)(X, Z)$, where $h = \frac{1}{2} \Ric _0 +
\frac1{24}\Scal\,g $ denotes the {\it normalized Ricci tensor} of $g$.

The normalized Ricci tensor is the form in which the Ricci tensor appears
in the well-known decomposition of the Riemannian curvature
$R = h \wedge \Id + W$, cf.~e.g. \cite{besse}; therefore, via the
(differential) Bianchi identity, the $\pm$-selfdual components, $C ^{\pm}$,
of the Cotton--York tensor are linked to the $\pm$-selfdual components of
the Weyl tensor by $ \delta W ^+ = C ^+$ and $ \delta W ^- = C ^-$.

Definition \ref{def1} can thus be rephrased as follows:
\begin{defi} \label{def2} A K{\"a}hler surface is {\rm weakly selfdual}
if its Cotton--York tensor is selfdual.
\end{defi}

The normalized Ricci tensor plays a natural role throughout this work.
For this reason, to simplify formulae, we write $s=\frac1{6} \Scal$
for the {\it normalized scalar curvature}, which is the trace of $h$.

\begin{lemma} \label{lemma1} {\rm (\cite{mats-tanno,jelonek})} For any
K{\"a}hler surface $(M, g,J,\omega)$ we have
\begin{equation} \label{eq:bianchirho}
\nabla _X \rho _0 = - 2 C ^- (JX) - \tfrac12 ds (X) \omega +
\tfrac12 \bigl( ds \wedge J X ^{\flat} - J ds \wedge X ^{\flat} \bigr).
\end{equation}

In particular, the K\"ahler surface $(M, g, J,\omega)$ is weakly selfdual if
and only if the following {\rm Matsumoto--Tanno identity}
\begin{equation} \label{eq:MT}
\nabla _X \rho _0 = - \tfrac{1}{2} ds (X) \omega
+ \tfrac{1}{2} ( ds \wedge J X ^{\flat} - J ds \wedge X ^{\flat})
\end{equation}
is satisfied \textup(for any vector field $X$\textup).
\end{lemma}
\begin{proof}
Since $2 h = \Ric - s\, g$, the Cotton--York tensor of any K{\"a}hler surface
can be written as follows:
\begin{equation*} \notag \begin{split}
2 C _{X, Y}(Z) &= - 2(\nabla _X h)(Y, Z) +  2(\nabla _Y h) (X, Z) \\
&= - (\nabla _X \Ric) (Y, Z) + (\nabla _Y \Ric ) (X, Z)
+  ds (X) \g{Y, Z} {-}  ds (Y) \g{X, Z} \\
&= - (\nabla _X \rho ) (Y, JZ) + (\nabla _Y \rho  ) (X, JZ)
+  ds (X) \g{Y, Z} -  ds (Y) \g{X, Z} \\
& =  (\nabla _{JZ} \rho ) (X, Y)  +  ds (X) \g{Y, Z} -  ds (Y) \g{X, Z}.
\end{split} \end{equation*}
(In order to obtain the last line from the preceding one, we use the fact
that $d \rho = 0$.) We then have:
\begin{equation*} 
(\nabla _Z \rho) (X, Y) = - 2 C  _{X, Y}(JZ)  - ds (X) \g{JY, Z}
+  ds (Y) \g{JX, Z},
\end{equation*}
or, equivalently
\begin{equation}\label{eq:nablarho0} \begin{split}
(\nabla _Z \rho_0) (X, Y)  = & - 2 C _{X, Y}(JZ) \\
& - \tfrac{3}{2} ds(Z)\g{JX,Y} - ds (X) \g{JY, Z} +  ds (Y) \g{JX, Z}.
\end{split}
\end{equation}
The anti-selfdual component of (\ref{eq:nablarho0}) gives identity
\eqref{eq:bianchirho}; the last statement follows immediately.
\end{proof}

\subsection{Twistor $2$-forms and K{\"a}hler metrics}

On any Riemannian manifold $(M,g)$, if $\Phi$ is an anti-selfdual $2$-form,
i.e., a section of the vector bundle $\Lambda ^- M$, then $\nabla\Phi$
is a section of the vector bundle $T ^* M \otimes \Lambda ^- M$.
This bundle has an orthogonal direct sum decomposition
\begin{equation} \label{eq:split} T ^* M \otimes \Lambda ^- M
= V ^{(0)} M \oplus V ^{(1)} M, \end{equation} in accordance with the
algebraic decomposition ${\R} ^4 \otimes \Lambda ^- {\R} ^4 = {\R} ^4
\oplus (\Sigma _+ \otimes \Sigma ^3 _-)$ into irreducible
sub-representations under the action of the orthogonal group. In
\eqref{eq:split}, $V ^{(0)} M$ corresponds to the factor ${\R} ^4$, hence
is isomorphic to $T ^* M$, whereas $V ^{(1)} M$ corresponds to the factor
$\Sigma _+ \otimes \Sigma ^3 _-$, so it is the kernel of the natural
contraction $T^* M\otimes\Lambda^- M\to T^* M$. The projection of the
connection to $V^{(0)} M$ may be identified with the exterior derivative or
divergence on anti-selfdual $2$-forms (which are related by the Hodge $*$
operator), while the projection to $V^{(1)} M$ is often called the {\it
twistor} or {\it Penrose operator} on anti-selfdual $2$-forms.

\begin{defi}\label{twistdef} An anti-selfdual $2$-form $\Phi$ is called a
{\rm twistor $2$-form} if $\nabla\Phi$ is a section of the sub-bundle
$V^{(0)}M$ of $T^*M \otimes \Lambda ^- M$.
\end{defi}

Any non-vanishing section $\Phi$ of $\Lambda ^- M$ can be written uniquely
as $\Phi = \lambda \, \omega _I$, where $\lambda = \frac{|\Phi|}{\sqrt{2}}$
is a positive function and $\omega _I$ is the K{\"a}hler form of an
anti-selfdual almost-complex structure $I$ on $(M, g)$.

\begin{lemma}\label{lem:twistor}\textup{\cite{pontecorvo}}
If $\Phi = \lambda \,\omega _I$ is a
non-vanishing section of $\Lambda ^- M$, then $\Phi$ is a twistor $2$-form
if and only if the almost-Hermitian pair $(\bar{g} = \lambda ^{-2} g, I)$
is K{\"a}hler, with K{\"a}hler form $\bar\omega=\lambda^{-2}\omega_I=
\lambda^{-3}\Phi$.
\end{lemma}
\begin{proof}
$\Phi$ is a twistor $2$-form if and only if there exists a $1$-form
$\gamma$ such that, at each point of $M$, $\nabla \Phi = \sum _{i = 1} ^3
I_i \gamma \otimes \omega _i$, where the triple $( I _1,I _2, I _3 = I _1 I
_2 )$ is any positively oriented, orthonormal frame of (anti-selfdual)
almost-complex structures at that point, and $\omega_i$ is the K{\"a}hler
form of $I_i$.

If $\Phi=\lambda\,\omega_I$ is a non-vanishing twistor $2$-form on $M$, then,
by choosing $I _1 = I$, we have
\begin{equation*}
\lambda \nabla I = (I \gamma - d \lambda) \otimes I +
I _2 \gamma \otimes I _2 + I _3 \gamma \otimes I _3.
\end{equation*}
Since the norm of $I$ is constant, this equality implies that $I \gamma = d
\lambda$. Now observe that the equation
\begin{equation}\label{eq:nablaI}
\lambda \nabla I = I _3 d \lambda \otimes I _2 - I _2 d\lambda \otimes I _3
\end{equation}
is equivalent to $I$ being parallel with respect to the Levi-Civita
connection of $\bar g=\lambda^{-2}g$. Hence if $\Phi$ is a twistor $2$-form
$(\bar g,I)$ is K{\"a}hler. Conversely, if $(\bar g,I)$ is K{\"a}hler then
\eqref{eq:nablaI} holds, from which it follows that $\nabla \Phi=\nabla
(\lambda\, \omega_I )$ is the section of $V ^{(0)} M$ corresponding to the
$1$-form $\gamma=-I d\lambda$.
\end{proof}

If $(M,g,J,\omega)$ is a K{\"a}hler surface then a $2$-form is anti-selfdual
if and only if it is trace-free and $J$-invariant, and there is the
following reformulation of Definition~\ref{twistdef}.

\begin{lemma}\label{lem:nablaphi0} Let $(M,g,J,\omega)$ be a K{\"a}hler
surface and let $\varphi_0$ be an anti-selfdual $2$-form . Then $\varphi_0$
is a twistor $2$-form if and only if there is a $1$-form $\beta$ such that
\begin{equation}\label{eq:nablaphi0}
\nabla_X\varphi_0= - \beta (X) \omega +
\beta \wedge J X ^{\flat} - J \beta \wedge X ^{\flat}
\end{equation}
for any vector field $X$.
\end{lemma}
\begin{proof} It is easily checked that the right hand side of
\eqref{eq:nablaphi0} is (the contraction with $X$ of) the general form of a
section of $V^{(0)}M\cong T^*M$.
\end{proof}

By the contracted Bianchi identity, $C^-$ is a section of $V^{(1)}M$, and
so the last statement of Lemma \ref{lemma1} can be rephrased as follows.

\begin{lemma} \label{lemma1bis}
A K{\"a}hler surface is weakly selfdual if and only if the trace-free part
$\rho _0$ of the Ricci form $\rho$ is a twistor $2$-form.
\end{lemma}
Together with Lemma~\ref{lem:twistor}, this implies:
\begin{prop} \label{proprho0} On the open set $M _0$ where $\rho _0 =\lambda
\,\omega_I$ does not vanish, a K{\"a}hler surface $(M,g,J,\omega)$ is weakly
selfdual if and only if the pair $(\bar{g} = \lambda ^{-2} g, I)$ is
K{\"a}hler.
\end{prop}

In particular it follows that on $M_0$, the selfdual Weyl tensor of $\bar
g$, with the orientation induced by $I$, is degenerate, and equal to
$\frac12 \bar s \,\bar\omega\otimes_0\bar\omega$, where $\bar s$ is the
(normalized) scalar curvature of $\bar g$,
$\bar\omega=\lambda^{-2}\omega_I$, and $\bar\omega\otimes_0\bar\omega$
stands for the traceless part of $\bar\omega\otimes\bar\omega$ viewed as an
endomorphism of $\Lambda^+ M$.  By the conformal covariance of the Weyl
tensor, it follows that on $M_0$, the \emph{anti-selfdual} Weyl tensor of
$g$ (using the orientation induced by $J$) is given by
\begin{equation}
W^- = \kappa\, \omega_I\otimes_0\omega_I
\end{equation}
where
\begin{equation} \label{kappabars}
\kappa = \bar{s} \lambda ^{-2}
\end{equation}
is the conformal scalar curvature of the Hermitian pair $(g, I)$, which is
related to the Riemannian scalar curvature of $g$ by
\begin{equation} \label{eq:kappas}
\kappa - s  = \delta \theta - |\theta| ^2;
\end{equation}
here $\theta$ denotes the {\it Lee form} of the pair $(g, I)$, defined by
$d \omega _I = - 2 \theta \wedge \omega _I$, see \cite{derdzinski,AG1}.
Note that we have normalized the conformal scalar curvature by a factor
$\frac16$ to be consistent with our normalization of the scalar curvature
$s$.

Notice that this only uses the fact that $M_0$ admits a non-vanishing
twistor $2$-form, namely $\rho_0$. On the other hand, $\rho_0$ is not an
arbitrary twistor $2$-form: by Lemma~\ref{lemma1}, the $1$-form $\beta$
defined by $\nabla\rho_0$ using~\eqref{eq:nablaphi0}, is equal to $\frac12
ds$, and so is exact. This fact, which is equivalent to the fact that the
Ricci form $\rho=\rho_0+ \frac32s\,\omega$ is closed, will be exploited in
the next subsection.

\subsection{Hamiltonian $2$-forms}

\begin{defi} A {\rm \biham/} $2$-form on a K{\"a}hler surface
$(M,g,J,\omega)$ is a closed $J$-invariant $2$-form $\varphi$ whose
trace-free \textup(i.e., anti-selfdual\textup) part $\varphi_0$ is a
twistor $2$-form.
\end{defi}

On a weakly selfdual K{\"a}hler surface the Ricci form $\rho$ is
\biham/. In general, \biham/ $2$-forms are characterized
by an analogue of the Matsumoto--Tanno identity~\eqref{eq:MT}.

\begin{lemma} \label{lem:genMT} Let $(M,g,J,\omega)$ be a K{\"a}hler surface.
Then a $J$-invariant $2$-form $\varphi=\varphi_0+\frac32\sigma\,\omega$ is
\biham/ if and only if
\begin{equation}\label{eq:genMT}
\nabla_X\varphi_0  = - \tfrac12 d\sigma (X) \omega +
\tfrac12 \bigl( d\sigma \wedge JX^{\flat} - J d\sigma \wedge X^{\flat} \bigr)
\end{equation}
for any vector field $X$.
\end{lemma}
\begin{proof}
This is immediate from Lemma~\ref{lem:nablaphi0}: in \eqref{eq:nablaphi0},
$d\varphi_0=-\frac32 d\sigma\wedge\omega$ if and only if $\beta=\frac12
d\sigma$.
\end{proof}

In order to explain the use of the term `\biham/', we recall the
following definition.

\begin{defi} A real function $f$ on a K{\"a}hler manifold $(M,g,J,\omega)$
is a \textup(real\textup) {\rm holomorphic potential} if the gradient ${\rm
grad} _g f$ is a holomorphic vector field, i.e., preserves $J$\textup;
equivalently, $f$ is a {\rm holomorphic potential} if $J {\rm grad} _g f$
is a Killing vector field with respect to $g$.

A holomorphic potential $f$ is therefore a momentum map for a hamiltonian
Killing vector field with respect to the symplectic form $\omega$.
\end{defi}

To any $J$-invariant $2$-form $\varphi=\varphi_0+\frac32 \sigma\,\omega$,
we may associate a {\it normalized} $2$-form
$\tilde\varphi=\frac12\varphi_0+\frac14 \sigma\,\omega$. For example, if
$\varphi$ is the Ricci form $\rho$, then $\tilde\varphi$ is the $2$-form
$\tilde\rho$ associated to the normalized Ricci tensor:
$\tilde\rho(\cdot,\cdot)=h(J\cdot,\cdot)$.

We are going to show that if $\varphi$ is \biham/, then the trace
and pfaffian of $\tilde\varphi$ are holomorphic potentials.

Recall that, in general, the {\it pfaffian} ${\rm pf} (\psi)$ of a $2$-form
$\psi$ is defined by $\frac12{\rm pf} (\psi) = {* (\psi \wedge \psi)}$,
where $*$ is the Hodge operator; alternatively, $\psi \wedge \psi = \frac14
{\rm pf} (\psi) \, \omega\wedge\omega$.

Since $\tilde\varphi=\frac12\varphi_0+\frac14\sigma\,\omega$, its pfaffian
$\pi$ is given by
\begin{equation} \label{eq:pfaffian}
\pi = \frac14 {\sigma^2} - \frac 12 |\varphi _0| ^2
= \Bigl(\frac{\sigma}{2} + \lambda\Bigr)\Bigl( \frac{\sigma}{2}
- \lambda\Bigr),
\end{equation}
where $\lambda=\frac{|\varphi_0|}{\sqrt2}$. We write this product as
$\pi=\xi\eta$ so that $\sigma=\xi+\eta$ and $\lambda=\frac12(\xi-\eta)$.

\begin{prop} \label{prop:biham} Let $(M,g,J,\omega)$ be a K{\"a}hler surface
and let $\varphi=\varphi_0+\frac32\sigma\,\omega$ be a \biham/
$2$-form.

Then the trace $\sigma$ and the pfaffian $\pi$ of
$\tilde\varphi=\frac12\varphi_0+\frac14\sigma\,\omega$ are Poisson-commuting
holomorphic potentials.
\end{prop}
\begin{proof} (i) Identity \eqref{eq:genMT} can be written in terms of
$\tilde\varphi$ as
\begin{equation} \label{eq:tildeMT}
\nabla _X \tilde\varphi = \tfrac{1}{4} \bigl( d\sigma \wedge JX^{\flat}
- J d\sigma \wedge X^{\flat} \bigr).
\end{equation}
Differentiating again and skew-symmetrizing, we get
\begin{equation*} \begin{split}
R _{X, Y}\cdot\tilde\varphi = [R _{X, Y}, \tilde\varphi] &= \tfrac{1}{4}
(\nabla _Y d\sigma \wedge JX^{\flat} - \nabla _X d\sigma \wedge JY^{\flat} \\
& \qquad - J \nabla _Y d\sigma \wedge X^{\flat}
+ J \nabla _X d\sigma \wedge Y^{\flat}).
\end{split}
\end{equation*}
Since $[R _{X,Y}, \tilde\varphi]$ is  $J$-invariant in $X, Y$, it follows that
\begin{multline*}
\nabla _Y d\sigma \wedge JX - J \nabla _Y d\sigma \wedge X
- \nabla _X d\sigma \wedge JY + J \nabla _X d\sigma \wedge Y = \\
- \nabla _{JY} d\sigma \wedge X - J \nabla _{JY} d\sigma \wedge JX
+ \nabla _{JX} d\sigma \wedge Y + J \nabla _{JX} d\sigma \wedge JY,
\end{multline*}
hence
\begin{equation} \label{eq:S=0}
S(Y) \wedge JX - J S (Y) \wedge X - S (X) \wedge JY + J S (X) \wedge Y = 0,
\end{equation}
where
\begin{equation*} S (X) = \nabla _X d\sigma + J \nabla _{JX} d\sigma.
\end{equation*}

As an algebraic object, $S$ is a symmetric, $J$-anti-commuting,
endomorphism of $TM$; hence by contracting \eqref{eq:S=0} with a vector field
$Z$ and taking the trace over $Y$ and $Z$, we see that $S=0$, and
therefore, $\sigma$ is a holomorphic potential.

(ii) From \eqref{eq:tildeMT} we derive
\begin{equation*}
(\nabla _X \tilde\varphi) \wedge \tilde\varphi
= \tfrac{1}{2} \bigl( X \wedge J d\sigma \wedge \tilde\varphi \bigr),
\end{equation*}
hence
\begin{equation} \label{eq:dpi}
d\pi = \tfrac12 {*} (J d\sigma \wedge \tilde\varphi).
\end{equation}
\noindent From this, we infer (again using \eqref{eq:tildeMT}):
\begin{equation} \label{eq:hessp} \begin{split}
\nabla_X d\pi &=
\tfrac12 \,{*} (J \nabla _X d\sigma \wedge \tilde\varphi
+  J d\sigma \wedge \nabla _X \tilde\varphi) \\
&= \tfrac12 \,{*}(J \nabla _X d\sigma \wedge \tilde\varphi
+ \tfrac{1}{4} J d\sigma \wedge d\sigma \wedge JX).
\end{split} \end{equation}
The second term of the right hand side of \eqref{eq:hessp} is clearly
$J$-invariant; the first term is $J$-invariant as well since $\sigma$ is a
holomorphic potential and $\tilde\varphi$ is $J$-invariant. Hence $\pi$ is
also a holomorphic potential.

(iii) By contracting equation~\eqref{eq:dpi} with $J d\sigma$, we see that
$\sigma$ and $\pi$ Poisson-commute.
\end{proof}

Since $\sigma$ and $\pi$ Poisson-commute, the Killing vector fields $K_1 :=
J {\rm grad} \,\sigma$ and $K_2 := J {\rm grad} \,\pi$ commute. Also
$\omega(K_1,K_2)=0$.

\begin{rem}\label{rem:hol} {\rm The Killing vector fields $K_1$ and $K_2$
need not be non-zero or independent in general. In particular $K_1$ and
$K_2$ both vanish if $\sigma$ is constant, since $\varphi$ is then parallel
by~\eqref{eq:genMT}.

Notice, however, that $K _1 ^{1,0}$, $K _2 ^{1,0}$, and hence $K _1 ^{1,0}
\wedge K _2 ^{1,0}$, are holomorphic, so that on each connected component
of $M$ there are three possibilities:
\begin{enumerate}
\item $K _1 \wedge K _2$ is non-vanishing on a dense open set.
\item $K _1$ is non-vanishing on a dense open set, but $K _1 \wedge K _2$
vanishes identically;
\item $K _1$ and $K _2$ vanish identically;
\end{enumerate}}
\end{rem}
One consequence of this remark is the following lemma.
\begin{lemma} \label{lem:M0dense} If $\varphi$ is a \biham/ $2$-form
on a K{\"a}hler surface $M$, then the open set $M_0$, where $\varphi_0$ is
non-zero, is empty or dense in each connected component of $M$.
\end{lemma}
\begin{proof}
On each component of $M$ where $K _1 = J\mathrm{grad}_g\sigma$ is non-zero,
the set $U$ where $d\sigma$ is non-vanishing is dense, hence
$\nabla\varphi_0$ is non-vanishing on $U$ and so the zero set of
$\varphi_0$ is dense in this connected component. On the other hand if
$K_1$ is identically zero on a component, then $\varphi_0$ is parallel on
that component, hence identically zero or everywhere non-zero.
\end{proof}

When $K_1$ and $K_2$ vanish identically, $\varphi$ does not contain much
information about the geometry of $M$ (it could be just a constant multiple
of $\omega$, even zero). In the other two cases, however, we shall obtain an
explicit classification of K{\"a}hler surfaces with a \biham/ $2$-form.
The keys to these classifications are Proposition~\ref{prop:biham} and the
following observation.

\begin{prop} \label{prop:dxi-deta} Let $\varphi$ be a \biham/
$2$-form on a K{\"a}hler surface $M$ and write $\sigma=\xi+\eta$
and $\pi=\xi\eta$ for the trace and pfaffian of $\tilde\varphi$.

Then on each connected component of $M$ where $\varphi_0$ is
not identically zero, $d\xi$ and $d\eta$ are orthogonal.
\end{prop}
\begin{proof} The contraction of \eqref{eq:genMT} with
$\varphi_0$ yields
\begin{align*}
\g{\nabla_X\varphi_0,\varphi_0}&=
\tfrac12\bigl(\varphi_0(d\sigma,JX)-\varphi_0(J d\sigma,X)\bigr)\\
&=-\varphi_0(J d\sigma,X)
\end{align*}
and hence
\begin{equation}\label{eq:swap}
2\lambda\, d\lambda = d(\lambda^2)
=\tfrac12 d(|\varphi_0|^2)=-\varphi_0(J d\sigma).
\end{equation}
Since $(\varphi_0\circ J)^2=\lambda^2\Id$, we deduce that $\varphi_0(J
d\lambda) =-\frac12\lambda\, d\sigma$ and therefore
\begin{equation*}\notag
\varphi _0\circ J \Bigl(\frac{d\sigma}{2} + d \lambda\Bigr) =
-\lambda \, \Bigl(\frac{d\sigma}{2} + d \lambda\Bigr), \quad
\varphi _0\circ J \Bigl(\frac{d\sigma}{2} - d \lambda\Bigr) =
\lambda \, \Bigl(\frac{d\sigma}{2} - d \lambda\Bigr).
\end{equation*}
This means that the 1-forms $d \xi$ and $d \eta$, wherever they are
non-zero, are eigenforms for the symmetric endomorphism $-\varphi _0\circ
J$, corresponding to the eigenvalues $\lambda$ and $-\lambda$,
respectively; in particular, they are orthogonal on the open set $M_0$
where $\lambda$ (i.e., $\varphi_0$) is non-zero. However, by
Lemma~\ref{lem:M0dense} this open set is empty or dense in each connected
component of $M$, and the result follows.
\end{proof}

\begin{rem}\label{rem:Idsig} {\rm
On the open set $M_0$ where $\lambda$ is non-zero, so that
$\varphi_0=\lambda\,\omega_I$, observe that equation~\eqref{eq:swap}
may be rewritten
\begin{equation}\label{eq:Idsig}
I d\sigma=2J d\lambda.
\end{equation}
Indeed, supposing only that $\varphi_0$ is a twistor $2$-form,
$\bar\omega=\lambda^{-3}\varphi_0$ is closed, and so, if
$\varphi=\varphi_0+\frac32\sigma\,\omega$,
\begin{align*}
d\varphi=d\bigl(\varphi_0+\tfrac32\sigma\,\omega\bigr)
&=3\lambda^2d\lambda\wedge\bar\omega+\tfrac32 d\sigma\wedge\omega\\
&=3\bigl(d\lambda\wedge\omega_I+\tfrac12 d\sigma\wedge\omega).
\end{align*}
Hence equation~\eqref{eq:Idsig} holds if and only if $\varphi$ is closed.}
\end{rem}

\subsection{Bi-extremal K{\"a}hler surfaces}

In the case that $(M,g,J,\omega)$ is weakly selfdual and $\varphi_0=\rho_0$,
we can take $\varphi=\rho$, so that $\sigma$ is the (normalized) scalar
curvature $s$, and $\pi$, the pfaffian $p$ of the normalized Ricci form
$\tilde\rho(\cdot, \cdot) = h (J \cdot, \cdot)$; evidently, $\tilde\rho =
\frac{1}{2} \rho _0 + \frac{1}{4}s \, \omega$.

Recall that a K{\"a}hler metric is said to be {\it extremal} if the
scalar curvature is a holomorphic potential.

\begin{defi}
A K{\"a}hler metric is called {\rm bi-extremal} if both the the
\textup(normalized\textup) scalar curvature
$s=\mathrm{tr}_\omega\tilde\rho$ and the pfaffian
$p=\mathrm{pf}(\tilde\rho)$ of the normalized Ricci form $\tilde\rho$ are
holomorphic potentials.
\end{defi}

Note that the potential function $p$ appearing in the above definition is
{\it not} the pfaffian of the {\it usual} Ricci form $\rho=\rho_0+\frac32
s\,\omega$; thus, our definition for bi-extremality differs from the one
given in \cite{maschler, hwang-maschler} (compare Theorem \ref{theo3} below
and \cite[Th.1.1 \& Prop.3.8]{hwang-maschler}).

Proposition~\ref{prop:biham} immediately implies:

\begin{prop}\label{prop:biextremal} A weakly selfdual K{\"a}hler metric is
bi-extremal.
\end{prop}

On a bi-extremal K{\"a}hler surface, the holomorphic potentials $s$ and $p$
automatically Poisson-commute, since $K_2$ preserves $g$, hence
$s$, so that $ds(K_2)=0$.

For a weakly selfdual K{\"a}hler surface, Proposition~\ref{prop:dxi-deta}
generically implies that $d\xi$ and $d\eta$ are orthogonal, where
$s=\xi+\eta$ and $p=\xi\eta$. We shall see in Sections~\ref{s:baks} and
\ref{s:coh1} that a bi-extremal K{\"a}hler surface satisfying this
orthogonality condition is weakly selfdual.

\subsection{The Bach tensor}\label{s:bach} The {\it Bach tensor}, $B$,  of an
$n$-dimensional Riemannian manifold $(M, g)$ is defined by
\begin{equation}\label{bach}
B _{X, Y}  =  \sum _{i = 1} ^n \big( - (\nabla _{e _i} C )_{e _i, X}(Y)
+ (W _{e _i, X} h (e _i), Y) \big),
\end{equation}
where, we recall, $C$ is the Cotton--York tensor and $h$ is the normalized
Ricci tensor (here, $\{ e _i \}$ is an arbitrary $g$-orthonormal frame).
When $n = 4$, the Bach tensor is conformal covariant of weight $-2$,
i.e., $B ^{\phi ^{-2} g} = \phi ^2 B ^g$, and $B$ can be indifferently
expressed in terms of $W ^+$ or of $W ^-$.  Specifically
\begin{equation} \label{bach+-} \begin{split}
B _{X, Y} &  = 2 \sum _{i = 1} ^n
\bigl( -(\nabla _{e _i} C ^+) _{e _i, X}(Y)
+ (W ^+  _{e _i, X} h (e _i), Y) \bigr) \\
& = 2 \sum _{i = 1} ^n \bigl( - (\nabla _{e _i} C ^-) _{e _i, X}(Y)
+ (W ^- _{e _i, X} h (e _i), Y)\bigr);
\end{split} \end{equation}
in particular, the Bach tensor vanishes whenever $W ^+$, $W^-$ or $\Ric _0$
vanishes.

If $(M,g,J)$ is a K{\"a}hler surface, the Bach tensor is easily computed by
using the above identity and the fact that $W ^+ = \frac{1}{2}
s\, \omega\otimes_{0} \omega$.  Indeed, if $B ^+$ and $B ^-$, denote the
$J$-invariant and $J$-anti-invariant parts of $B$, we get
\begin{equation}
B ^+ = s \, \Ric_0 + 2 (\nabla ds)^+_0, \qquad
B ^- = - (\nabla d s)^-,
\end{equation}
where $(\nabla ds)^+_0$ is the $J$-invariant trace-free part of the Hessian
and $(\nabla d s)^-$ is the $J$-anti-invariant part.  (This formula
is due to Derdzi\'nski \cite{derdzinski}, where it was obtained by a
different argument.)

It follows that $B$ is $J$-invariant if and only if $(M,g,J)$ is an
extremal K{\"a}hler surface; moreover, if this holds, we get $B= s\, \Ric_0
+ 2 (\nabla ds)_0$.  On the open set $U$ where $s$ has no zero, the
vanishing of $B $ then means that the conformally related metric $\tilde{g}
=s^{-2}g$ is Einstein. Therefore, on $U$, the following two statements are
thus equivalent (see also \cite{derdzinski}):
\begin{enumerate}
\item The Bach tensor of the K{\"a}hler surface $(M, g, J)$ vanishes;
\item $(M, g, J)$ is extremal and the conformally related metric
$\tilde{g} =s^{-2}g$ is Einstein.
\end{enumerate}
Note also that when $B$ is $J$-invariant, it is determined by the
associated anti-selfdual $2$-form ${\tilde B}(\cdot,\cdot) =
B(J\cdot,\cdot)$, which is also given by
\begin{equation}\label{eq:Bachform}
{\tilde B}=(d J d s)_0 +s \rho_0.
\end{equation}
On a weakly selfdual K{\"a}hler surface, $C^-=0$, while $W^-=\frac12
\kappa\,\omega_I\otimes_0\omega_I$. It follows that ${\tilde B}$ is a
multiple of $\kappa\rho_0$, which vanishes if and only if $W^-=0$ or
$\rho_0=0$
\begin{prop} A weakly selfdual K{\"a}hler surface is Bach-flat \textup(i.e.,
has vanishing Bach tensor\textup) if and only if it is selfdual or
K{\"a}hler--Einstein.
\end{prop}

\subsection{Rough classification of weakly selfdual K{\"a}hler surfaces}

We have seen in Lemma~\ref{lem:M0dense} that the open set $M_0$, on which
the trace-free part $\varphi_0$ of a \biham/ $2$-form is non-zero,
is empty or dense in each connected component of $M$. For a weakly
selfdual K{\"a}hler surface $(M,g,J,\omega)$, the Ricci form $\rho$ is
\biham/, and $M_0$ is the set of points at which $g$ is not
K{\"a}hler--Einstein. Because $\rho$ is closely linked to the anti-selfdual
Weyl tensor of $M$, we can obtain more information about $M_0$ except when
$g$ is selfdual.

We first recall the general fact, first observed by A. Derdzi\'nski in
\cite{derdzinski}, that for any K{\"a}hler surface $(M, g, J)$ with
non-vanishing scalar curvature $s$, the conformally related metric
$\tilde{g} = s ^{-2} g$ satisfies $\delta ^{\tilde{g}} W ^+ = 0$; moreover,
up to rescaling, $\tilde{g}$ is the unique metric in the conformal class
$[g]$ that satisfies this property. This follows from the fact that the
selfdual Cotton--York tensor $C ^+$ of any K{\"a}hler surface can be
written as
\begin{equation} \label{eq:C+s}
C ^+ (X) = - W ^+ \Bigl(\frac{ds}{s} \wedge X\Bigr);
\end{equation}
on the other hand, the selfdual Cotton--York tensors of two conformally
related metrics $g$ and $f ^{-2} g$ are related by
\begin{equation*} \label{eq:CY}
C ^{+,f ^{-2} g} (X)  = C ^{+,g} (X) + W ^+ \Bigl(\frac{d f}{f}\wedge X\Bigr);
\end{equation*}
it follows from \eqref{eq:C+s} that the selfdual Cotton--York tensor of the
metric $s ^{-2} g$ vanishes identically; moreover, as $W ^+$ has no kernel
(for $s$ non-vanishing), the latter property characterizes $s ^{-2} g$ up to
a constant multiple.

\begin{lemma} \label{propbars} On each connected component of the open set
$M _0$ where $\rho _0$ does not vanish, the scalar curvature $\bar{s}$ of
$\bar{g}$ is a constant multiple of $\lambda ^{-1}$, i.e.,
\begin{equation} \label{eq:C}
\bar{s} = c \lambda ^{-1},
\end{equation}
where $\lambda$ is the positive eigenvalue of $\Ric _0$ and $c$ is a
constant.
\end{lemma}
\begin{proof} We apply the preceding argument to the K{\"a}hler pair
$(\bar{g}, I)$ on $M_0$ and to $g = \lambda ^2 \bar{g}$; by hypothesis, $g$
satisfies $\delta^g W ^- = C ^- = 0$, where $W ^-$ is actually the {\it
selfdual} Weyl tensor of $g$ for the orientation induced by $I$; from the
above mentioned uniqueness property, it follows that, wherever $\bar s$ is
non-zero, $g$ coincides with $\bar{s} ^{-2} \bar{g}$ up to rescaling, i.e.,
that $\bar{s}$ is a locally constant multiple of $\lambda ^{-1}$.  However,
the same holds on the interior of the zero set of $\bar s$. Hence by the
continuity of $\bar s$ on $M_0$, $\bar s=c\lambda^{-1}$ for some constant
$c$ on each connected component of $M_0$.
\end{proof}

A more global statement may be obtained using the conformal scalar
curvature $\kappa=\bar{s} \lambda ^{-2}$. Since the anti-selfdual Weyl
tensor of $g$ is given by $W^-=\frac12\kappa\, \omega_I\otimes_0\omega_I$, it
follows that $\kappa^2$ is equal to $|W^-|^2$ on $M_0$, up to a numerical
factor, and hence we may extend $\kappa$ continuously to the closure of
$M_0$. Also $\lambda$ is globally defined and continuous.

Therefore, using the fact that the closure of $M_0$ is a union of connected
components by Lemma~\ref{lem:M0dense}, we can rewrite Lemma~\ref{propbars}.

\begin{lemma} \label{propkappa} Let $(M,g,J,\omega)$ be a weakly selfdual
K{\"a}hler surface. Then, on each component of $M$ where $\rho_0$ is not
identically zero, the conformal scalar curvature $\kappa$ of $(g, I)$ is
linked to $\lambda$ by
\begin{equation} \label{eq:Ckappa} \kappa\lambda^{3} = c,
\end{equation}
where $c$ is the constant of Lemma~{\rm \ref{propbars}}.
Moreover, $c=0$ if and
only if $W^-=0$ on that component.
\end{lemma}

This lemma yields the following rough classification of weakly selfdual
K{\"a}hler surfaces (see also \cite{jelonek}).

\begin{prop} \label{prop:tri} Let $(M, g, J,\omega)$ be a weakly selfdual,
connected, K{\"a}hler surface. Then either:

\begin{enumerate}
\item $\rho_0$ is identically zero so $(g,J)$ is K{\"a}hler--Einstein; or

\item the scalar curvature $s$ of $g$ is constant, but $\rho_0$ is
not identically zero; then, $(g, J)$ is locally the K{\"a}hler product of two
Riemann surfaces of constant curvatures; or

\item $s$ is not constant and $g$ is selfdual; or

\item $W ^-$ and $\rho _0$ have no zero: then, the K{\"a}hler
metric $(\bar{g} = \lambda ^{-2} g, I)$ of Proposition {\rm \ref{proprho0}}
is extremal and globally defined on $M$; in particular, $W ^-$ is
degenerate everywhere.
\end{enumerate}
\end{prop}
\begin{proof} If $s$ is constant, then by \eqref{eq:MT} the Ricci form is
parallel. Hence either $g$ is locally irreducible, and is Einstein, or $(g,
J)$ is locally the K{\"a}hler product of two Riemann surfaces of constant
curvatures.

If $s$ is not constant, then by Lemma~\ref{lem:M0dense} the open set $M _0$
where $\rho _0$ does not vanish is an open dense subset of $M$.  However,
by Lemma~\ref{propkappa}, $\kappa\lambda^3$ is constant. If this constant
is zero, then $\kappa$ must vanish identically and $M$ is selfdual;
otherwise $\kappa$ and $\lambda$ have no zero on $M$, so $M _0 = M$, the
K{\"a}hler pair $(\bar{g} = \lambda ^{-2} g, I)$ is defined on $M$ and $W
^-$ is degenerate, but nonzero everywhere. As observed in Section 1.5,
for a weakly selfdual \ka surface the Bach form  ${\tilde B}$ is a
multiple of $\omega_I$. Since $B$ is a conformally covariant tensor, it
follows that Bach tensor of ${\bar g}$ is $I$-invariant, showing that
$({\bar g},I)$ is an extremal \ka metric.
\end{proof}

K{\"a}hler--Einstein metrics and K{\"a}hler products of Riemann surfaces
clearly \emph{are} weakly selfdual. Since these are well studied, we
henceforth assume that $s$ is not constant (on any component), i.e., $K _1$
is non-vanishing on a dense open set.  In Section~\ref{s:baks}, we analyse
the generic case that $K _1$ and $K _2$ are independent, while
Section~\ref{s:coh1} is devoted to the case that $K_1\wedge K_2$ vanishes
identically (but $K_1$ is non-zero). In both sections, we obtain an explicit
local classification within the more general framework of K{\"a}hler
surfaces with a \biham/ $2$-form.

\section{Ortho-toric K{\"a}hler surfaces}\label{s:baks}

\subsection{Toric K{\"a}hler surfaces}

We have seen in Section \ref{s:wsdk} that on a K{\"a}hler surface with a
\biham/ $2$-form $\varphi$---in particular on a weakly selfdual
K{\"a}hler surface---the trace $\sigma$ and the pfaffian $\pi$ of the
associated normalized $2$-form $\tilde\varphi$ are holomorphic potentials
for hamiltonian Killing vector fields $K_1 = J {\rm grad} \, \sigma$ and
$K_2 = J {\rm grad} \, \pi$. Furthermore, $\sigma$ and $\pi$
Poisson-commute, i.e., $\omega(K_1,K_2)=0$.

A (usually compact) K{\"a}hler surface $(M,g,J,\omega)$, with holomorphic
Killing vector fields $K_1$ and $K_2$ which are independent on a dense open
set and satisfy $\omega(K_1,K_2)=0$, is said to be {\it toric}.  We begin
this section by recalling the local theory of such surfaces, and we
therefore assume that $K_1$ and $K_2$ are everywhere independent and that
$x_1$ and $x_2$ are globally defined momentum maps for $K_1$ and $K_2$.

The condition $\omega(K_1, K_2)=0$ is equivalent to the fact that $x_1$ and
$x_2$ commute for the Poisson bracket determined by $\omega$. Hence also
$[K_1, K_2]=0$, and since $K_1$, $K_2$, $J K_1$ and $J K_2$ are all
holomorphic, they all commute. In particular the rank $2$ distributions
$\Pi$, generated by $K_1$ and $K_2$, and $J\Pi$, generated by $J K_1$ and
$J K_2$, are integrable.

These distributions $\Pi, J\Pi$ are also orthogonal, since $\g{J K_1,
K_2}=0$. It follows that $K_1$, $K_2$, $J K_1$ and $J K_2$ form a
frame. Since they commute, the $1$-forms in the dual coframe are closed,
and may be written $d t_1, d t_2, J dt_1, J dt_2$, where $t_1$ and $t_2$
are only given locally and up to an additive constant.  Now observe that
\begin{equation*}\begin{split}
J dt_1 &= \frac{|K_2|^2 dx_1-\g{K_1,K_2} dx_2}{|K_1\wedge K_2|^2}\\
J dt_2 &= \frac{|K_1|^2 dx_2-\g{K_1,K_2} dx_1}{|K_1\wedge K_2|^2}
\end{split}\end{equation*}
and so
\begin{equation*}
J dt_i = \sum_{j=1,2} G_{ij} dx_j\qquad(i=1,2),
\end{equation*}
where $G_{ij}$ is a positive definite symmetric matrix of functions of
$x_1$ and $x_2$ (note that $K_i=\partial/\partial t_i$).  These $1$-forms
are closed if and only if $G_{ij}$ is the Hessian of a function of $x_1$
and $x_2$. The following well-known explicit classification is then readily
obtained (see \cite{G:kstv}, \cite{Abreu}).

\begin{prop}\label{p:tak} Let $G_{ij}$ be a positive definite $2\times 2$
symmetric matrix of functions of $2$-variables $x_1,x_2$ with inverse
$G^{ij}$.  Then the metric
\begin{equation*}
\sum_{i,j} \bigl( G_{ij} dx_i dx_j + G^{ij} dt_i dt_j \bigr)
\end{equation*}
is almost-K{\"a}hler with K{\"a}hler form
\begin{equation*}
\omega = dx_1\wedge dt_1+dx_2\wedge dt_2
\end{equation*}
and has independent hamiltonian Killing vector fields $\partial/\partial
t_1, \partial/\partial t_2$ with Poisson-commuting momentum maps $x_1$ and
$x_2$. Any almost K{\"a}hler structure with such a pair of Killing vector
fields is of this form \textup(where the $t_i$ are locally defined up to an
additive constant\textup), and is K{\"a}hler if and only if $G_{ij}$ is the
Hessian of a function of $x_1$ and $x_2$.
\end{prop}

\subsection{The \biax/ case}

Propositions~\ref{prop:biham} and~\ref{prop:dxi-deta} motivate the
following definition.

\begin{defi} A K{\"a}hler surface $(M,g,J,\omega)$ is {\rm \biax/} if it
admits two independent hamiltonian Killing vector fields with Poisson-commuting
momentum maps $\xi\eta$ and $\xi+\eta$ such that $d\xi$ and $d\eta$ are
orthogonal.
\end{defi}

An explicit classification of \biax/ K{\"a}hler metrics follows from
Proposition~\ref{p:tak} by changing variables and imposing the
orthogonality of $d\xi$ and $d\eta$. However, since the coordinate change
is awkward, and we have not spelt out the proof of Proposition~\ref{p:tak},
we give a self-contained proof of this classification.

\begin{prop}\label{prop:biax}
The almost-Hermitian structure $(g,J,\omega)$ defined by
\begin{gather} \label{eq:newg} \begin{split}
g &= (\xi - \eta) \, \Bigl ( \frac{d \xi ^2}{F (\xi)}
- \frac{d \eta ^2}{G (\eta)} \Bigr)\\
&\qquad + \frac{1}{\xi - \eta} \bigl ( F (\xi) (dt + \eta dz) ^2
- G (\eta) (dt + \xi dz) ^2 \bigr),
\end{split}\\
\label{eq:newJ} \begin{split}
J d \xi &= \frac{F (\xi)}{\xi - \eta} (dt + \eta\, dz), \qquad
J dt = -\frac{\xi \,d \xi}{F (\xi)} - \frac{\eta \,d \eta}{G (\eta)}, \\
J d \eta &= \frac{G (\eta)}{\eta - \xi} (dt + \xi\, dz), \qquad
J dz = \frac{d \xi}{F (\xi)} + \frac{d \eta}{G (\eta)},
\end{split}\\
\label{eq:newomega}
\omega=d\xi\wedge(dt+\eta\, dz)+d\eta\wedge (dt +\xi\, dz)
\end{gather}
is an \biax/ K{\"a}hler structure for any functions $F,G$ of one variable.
Every \biax/ K{\"a}hler surface is of this form, where $t,z$ are locally
defined up to an additive constant.
\end{prop}
\begin{proof} (i) The K{\"a}hler form may be written
\begin{equation*}
\omega=d(\xi+\eta)\wedge dt +d(\xi\eta)\wedge dz
\end{equation*}
which is certainly closed. If is also immediate that $\partial/\partial t$
and $\partial/\partial z$ are hamiltonian Killing vector fields with
Poisson-commuting momentum maps $\xi+\eta$ and $\xi\eta$.  Since $dt+i Jdt$
and $dz+ i J dz$ are closed, $J$ is integrable, and the K{\"a}hler surface is
clearly \biax/.

(ii) Conversely, suppose that $(g,J,\omega)$ is an \biax/ K{\"a}hler surface
with Killing vector fields $K_1$, $K_2$. Since the dual frame to $K_1, K_2,
JK_1, JK_2$ consists of closed $1$-forms, we may write it as $dt, dz, J dt,
J dz$, where $t$ and $z$ are locally defined up to an additive constant.
Note also that $d\xi, d\eta,dt, dz$ are linearly independent
$1$-forms---where $\xi+\eta$ and $\xi\eta$ are the momentum maps of $K_1$
and $K_2$---so we may use $(\xi,\eta,t,z)$ as a coordinate system.

Since $(J dz) (K_1)=0$ and $(J dz) (K_2)=0$ we may write
\begin{equation*}
J dz = \frac{d\xi}{F}+\frac{d\eta}{G}
\end{equation*}
for some functions $F$ and $G$ (of $\xi$ and $\eta$).
The equations
\begin{align*}
0&=(J dz)(J K_1)=-\g{J dz, d\xi+d\eta}\\
\tag*{and} 1&=(J dz)(J K_2)=-\g{J dz, \eta \,d\xi+\xi \,d\eta}
\end{align*}
give $F=|d\xi|^2(\xi-\eta)$ and $G=|d\eta|^2(\eta-\xi)$, using the fact
that $d\xi$ and $d\eta$ are orthogonal. A similar argument tells us that
\begin{equation*}
J dt = \frac{\xi\,d\xi}{F}+\frac{\eta\,d\eta}{G}.
\end{equation*}
for the same functions $F$ and $G$. Now since $J dt$ and $J dz$ are closed,
we obtain $(\xi-\eta)F_\eta=0$ and $(\eta-\xi)G_\xi=0$, so that $F=F(\xi)$
and $G=G(\eta)$. The K{\"a}hler form $\omega$ is evidently given
by~\eqref{eq:newomega}, and since we know $Jdt$ and $Jdz$, we
readily obtain~\eqref{eq:newJ}, and hence the metric~\eqref{eq:newg}.
\end{proof}

Any \biax/ K{\"a}hler surface $(M, g, J,\omega)$ comes equipped with an
anti-selfdual almost-complex structure, $I$, whose K{\"a}hler form,
$\omega _I$ is defined by
\begin{equation} \label{eq:omegaI} \begin{split}
\omega _I & = \frac{d \xi \wedge J d \xi}{|d \xi| ^2}
- \frac{d \eta \wedge J d \eta}{|d \eta| ^2} \\
&= d \xi \wedge (dt + \eta dz) - d \eta \wedge (dt + \xi dz);
\end{split}
\end{equation}
equivalently
\begin{equation} \label{eq:newI} \begin{split}
I d \xi &= \hphantom{-} J d \xi =
\hphantom{-} \frac{F (\xi)}{\xi - \eta} (dt + \eta dz),\qquad
I dt = -\frac{\xi d \xi}{F (\xi)} + \frac{\eta d \eta}{G (\eta)}, \\
I d \eta &= - J d \eta = - \frac{G(\eta)}{\eta-\xi} (dt + \xi dz),\qquad
I dz = \frac{d \xi}{F (\xi)} - \frac{d \eta}{G (\eta)}.
\end{split} \end{equation}

\begin{prop} \label{prop:newIkaehler} For any \biax/ K{\"a}hler surface,
the almost-Hermitian pair $(\bar{g} = (\xi - \eta) ^ {-2} g, I)$ is
K{\"a}hler.
\end{prop}
\begin{proof} Clearly $I dt$ and $I dz$ are closed, so $I$ is integrable. From
\eqref{eq:omegaI}, we easily infer that the Lee form $\theta$ of the
Hermitian pair $(g, I)$, defined by $d\omega_I=-2\theta\wedge\omega_I$, is
\begin{equation} \label{eq:theta}
\theta = - d \log{|\xi - \eta|}.
\end{equation}
It follows that $\bar{\omega} := (\xi - \eta) ^ {-2} \omega _I$ is closed,
i.e., the pair $(\bar{g} = (\xi - \eta) ^ {-2} g, I)$ is K{\"a}hler.
\end{proof}

In particular, on any \biax/ K{\"a}hler surface, the anti-selfdual Weyl
tensor---which is the {\it selfdual} Weyl tensor of $g$ for the orientation
induced by $I$---is degenerate: $W^-=\kappa\,\omega_I\otimes_0\omega_I$,
where $\kappa$ is the conformal scalar curvature of the Hermitian pair
$(g,I)$.

\begin{rem} {\rm The vector fields $K _1$ and $K _2$ are still Killing
with respect to $\bar{g}$ and hamiltonian with respect to $\bar{\omega}$,
with momentum maps $- \frac{1}{\xi - \eta}$ and $ -
\frac{\xi + \eta}{2(\xi - \eta)}$ respectively. However, the K{\"a}hler
metric $({\bar g}, I)$ is {\it not} \biax/ in general, as it can be
checked using Lemma \ref{lem:curv} below. }
\end{rem}

Combining Propositions~\ref{prop:biham},~\ref{prop:dxi-deta}
with~\ref{prop:biax} and~\ref{prop:newIkaehler}, we obtain the following
theorem.

\begin{theo} A K{\"a}hler surface is \biax/ if and only if it admits
a \biham/ $2$-form whose associated Killing vector fields are
independent. The K{\"a}hler structure is then given explicitly in terms
of two arbitrary functions $F,G$ of one variable
by~\eqref{eq:newg}--\eqref{eq:newomega}.
\end{theo}

Indeed, using~\eqref{eq:omegaI}--\eqref{eq:newI}, notice that, by
definition, $Jd(\xi+\eta)=Id(\xi-\eta)$, cf.\ Remark~\ref{rem:Idsig},
so the \biham/ $2$-form $\varphi$ is
$\frac12(\xi-\eta)\omega_I+(\xi+\eta)\omega$.

\subsection{Ortho-toric weakly selfdual K{\"a}hler surfaces}

The curvature of an \biax/ K{\"a}hler surface is entirely determined by the
scalar curvature $s$ of $g$, the conformal scalar curvature $\kappa$ of the
Hermitian pair $(g,I)$, and the trace-free part $\rho _0$ of the Ricci form
of $(g, J)$.

\begin{lemma} \label{lem:curv} For any \biax/ K{\"a}hler surface
$(M, g, J,\omega)$, $\rho_0$ is a multiple $\mu$ of the K{\"a}hler form
$\omega_I$ of the Hermitian pair $(g, I)$, and $\mu,s,\kappa$ are
given by
\begin{align}
\label{eq:newmu}
\mu &=  \frac{F' (\xi) - G' (\eta)}{2(\xi - \eta) ^2}
- \frac{F'' (\xi) + G'' (\eta)}{4(\xi - \eta)},\\
\label{eq:newscalar}
s &= - \frac{F''(\xi) - G'' (\eta)}{6(\xi - \eta)},\\
\label{eq:newkappa}
\kappa &= - \frac{F'' (\xi) - G'' (\eta)}{6(\xi - \eta)}
+ \frac{F' (\xi) + G' (\eta)}{(\xi - \eta) ^2}
- \frac{2(F (\xi) - G (\eta))}{(\xi - \eta) ^3}.
\end{align}
In particular, on the open subset of $M$ where $\mu$ has no zero, the
anti-selfdual almost-complex structure determined by $\rho _0$ is equal to
$I$.
\end{lemma}
\begin{proof} From \eqref{eq:newomega}, we infer that the volume-form
$v _g = \frac{1}{2} \omega \wedge \omega$ of $g$ is given by
\begin{equation*}
v _g = -  (\xi - \eta) \, d \xi \wedge d \eta \wedge dt \wedge dz,
\end{equation*}
since $t$ and $z$ are the real parts of $J$-holomorphic coordinates.
By putting $v _0 = dt \wedge J dt \wedge dz \wedge J dz$, we have
\begin{equation*} \rho = - \frac{1}{2} d J d \log{\frac{v _g}{v _0}}.
\end{equation*}
Now according to \eqref{eq:newJ},
$$ v _0 = - \frac{\xi - \eta}{ F (\xi) G (\eta)} d \xi \wedge d \eta
\wedge dt  \wedge dz, $$
and hence
\begin{equation*} \frac{v _g}{v _0} =  F (\xi) G (\eta);
\end{equation*}
this implies
\begin{equation*}
\rho = -\frac{1}{2} d J d \log{|F(\xi)|} - \frac{1}{2} d J d \log{|G(\eta)|},
\end{equation*}
from which \eqref{eq:newmu} and \eqref{eq:newscalar} follow easily.

\par From \eqref{eq:kappas} and \eqref{eq:theta}, we get
\begin{equation*} \label{kappaint}
\kappa = s - 2 \, \frac{|d \xi| ^2 + |d \eta| ^2}{(\xi - \eta) ^2} -
\frac{\Delta (\xi - \eta)}{\xi - \eta};
\end{equation*}
on the other hand, we compute that
\begin{equation*} \label{deltaxi}
\Delta \xi = - \frac{F' (\xi)}{\xi - \eta}, \ \ \
\Delta \eta = \frac{G' (\eta)}{\xi - \eta};
\end{equation*} and we obtain \eqref{eq:newkappa}.
\end{proof}

\begin{prop} \label{prop:newextremal} An \biax/ K{\"a}hler surface $M$
is extremal if and only if $F$ and $G$ are of the form
\begin{equation} \label{eq:extremal} \begin{split}
F (x) &= k x ^4 + \ell x ^3 + Ax ^2 + B _1 x + C _1, \\
G (x) &= k x ^4 + \ell x ^3 + Ax ^2 + B _2 x + C _2,
\end{split} \end{equation}
in which case
\begin{equation}\label{eq:exts}
s =  - 2 k (\xi + \eta) -  \ell,
\end{equation}
and $({\bar g} = (\xi - \eta)^{-2}g, I)$ is an extremal \ka metric as
well.

Moreover, $M$ is
\begin{itemize}
\item Bach-flat if and only if $4k (C _1 - C _2) = (B _1 - B _2)\ell$;
\item of constant scalar curvature if and only if $k=0$;
\item scalar-flat \textup(i.e., anti-selfdual\textup) if
and only if $k=\ell=0$.
\end{itemize}
\end{prop}
\begin{proof} Since the scalar curvature $s$ is a function of $\xi$ and
$\eta$, $J {\rm grad} _g s$ belongs to the span of the Killing vector
fields $K _1$ and $K _2$ and commutes with them; if it is itself a Killing
vector field, it has to be a linear combination of $K _1$ and $K _2$ with
constant coefficients, i.e. $s = a (\xi + \eta) + b \xi \eta + c$, where
$a, b, c$ are constants. By \eqref{eq:newscalar}, this implies
\eqref{eq:extremal}.  Finally, using~\eqref{eq:Bachform}, we easily compute
that the anti-selfdual $2$-form associated to the Bach tensor of an
\biax/ extremal K{\"a}hler surface is
\begin{equation*}
{\tilde B}= \frac{4k (C _1 - C _2) - (B _1 - B _2)\ell}{2 (\xi - \eta)^2}
\,\omega _I.
\end{equation*}
Since $B$ is also $I$-invariant, the \ka metric  $({\bar g}, I)$
is extremal as well (see Section 1.5).
\end{proof}

\smallskip\noindent {\bf Example 1.}
For $k\neq0$ and $4(C_1-C_2)=(B _1 - B _2)\frac\ell k$, we obtain
explicit Bach-flat K{\"a}hler surfaces with non-constant scalar
curvature. These metrics are therefore not anti-selfdual, and for $B_1\neq
B_2$ they are not selfdual either (note that a selfdual K{\"a}hler surface is
bi-extremal and see the next Proposition). According to
Section~\ref{s:bach}, the metric
\begin{equation*}
{\tilde g} = (2 k (\xi + \eta) +  \ell)^{-2} g,
\end{equation*}
which is defined on the open subset where $2 k (\xi + \eta) +  \ell \neq 0$,
is Einstein, Hermitian (but non-K{\"a}hler) with a locally defined toric
isometric action.

\begin{prop} \label{prop:newbiextremal} An \biax/ K{\"a}hler surface $M$
is bi-extremal if and only if $F$ and $G$ are of the form
\begin{equation} \label{biextremal} \begin{split}
F (x) &= k x ^4 + \ell x ^3 + Ax ^2 + B x + C _1, \\
G (x) &= k x ^4 + \ell x ^3 + Ax ^2 + B x + C _2,
\end{split} \end{equation}
in which case the Ricci form is given by
$\rho=-2k\varphi-\ell\omega$. Hence $M$ is weakly selfdual and is
\begin{itemize}
\item selfdual if and only if $C_1=C_2$;
\item K{\"a}hler--Einstein if and only if $k=0$;
\item Ricci-flat if and only if $k=\ell=0$.
\end{itemize}
\end{prop}
\begin{proof} Since a bi-extremal surface is extremal, we may apply
Proposition~\ref{prop:newextremal}. We then compute
\begin{align}\label{eq:extmu}
\mu &= - k (\xi - \eta) + \frac{B _1 - B _2}{2 (\xi - \eta) ^2},\\
p  & = 4 k ^2 \xi \eta + k \ell (\xi + \eta) + \frac{\ell ^2}{4}
- k\, \frac{B _1 - B _2}{\xi-\eta} +\frac{(B _1 - B _2)^2}{4(\xi-\eta) ^4}.
\end{align}
As in the proof of Proposition~\ref{prop:newextremal}, $p$ cannot be a
$J$-holomorphic potential unless it is a linear combination of $\xi
+ \eta$ and $\xi \eta$ with constant coefficients; this in turn is
equivalent to the condition $B _1 = B _2$. Since $\mu=-2k\lambda$
and $s=-2k\sigma-\ell$ it follows that $\rho= -2k\varphi-\ell\omega$.
Hence $\rho$ is a \biham/ and so $M$ is weakly selfdual
by~\eqref{lemma1bis}.

By substituting in the expression of $\kappa$ given by \eqref{eq:newkappa}
we get $\kappa = - \frac{2 (C _1 - C _2)}{(\xi - \eta) ^3}$; since the
condition $W ^- = 0$ is equivalent to $\kappa = 0$, the characterization of
the selfdual case follows. Also $M$ is Einstein if and only if $\mu=0$, and
so the last two assertions are immediate.
\end{proof}

\begin{rem} {\rm For $k\neq 0$, we can set $k=-\frac12$ and $\ell=0$ by
a simultaneous affine change of $\xi$ and $\eta$. In this case
$\rho=\varphi$, and so the \biax/ reduction is defined by
$\rho$. However, not all weakly selfdual K{\"a}hler surfaces can be put in
\biax/ form; for example weakly selfdual metrics belonging to the general
family of cohomogeneity-one extremal metrics considered by
E. Calabi~\cite{calabi1} are not in general \biax/, since $K_1$ and $K_2$
are then collinear. We discuss this case in Section~\ref{s:coh1}.

On the other hand, the examples with $k=0$ in the above Proposition show
that among K{\"a}hler--Einstein metrics (which are weakly selfdual), there
are some which can be put into \biax/ form, because even though $\rho$ does
not define an \biax/ reduction, there happens to be another \biham/ $2$-form
$\varphi$. If additionally $A=\ell=0$ or $C_1=C_2$, these examples in fact
have cohomogeneity one and their \biax/ form arises from a choice of maximal
torus inside the isometry group.  However, the other examples do not have
additional symmetries and therefore appear to be new.}
\end{rem}

\begin{prop}\label{propnewsymmetries} On an \biax/ extremal K{\"a}hler
surface $M$, the space of infinitesimal symmetries of the K{\"a}hler
structure is generated by $K _1$ and $K _2$ \textup(and infinitesimal
rotations in this plane if they are globally defined\textup), except in the
following two cases:
\begin{enumerate}
\item $M$ is locally a complex space form, i.e., K{\"a}hler--Einstein and
selfdual; or

\item $M$ is Ricci-flat \textup(hence anti-selfdual\textup) and,
locally, of cohomogeneity one.

In terms of $F$ and $G$, these two cases are respectively described by
\begin{equation} \label{eq:cond1}
F(x) = G(x) = \ell x^3 + Ax^2 + Bx + C;
\end{equation}
and
\begin{equation} \label{eq:cond2} \begin{split}
&F(x) =  Bx + C_1, \\ &G(x) = Bx + C_2.
\end{split} \end{equation}
\end{enumerate}
\end{prop}
\begin{proof} Suppose there exist a third infinitesimal symmetry of
$M$, say $K _3$, which does not lie in the plane spanned by $K _1$ and $K
_2$; then, we must have $ds(K_i) = 0$ and $ d \mu (K_i) = 0$, for $i = 1,
2, 3$; this implies that $ds$ and $d \mu$ are colinear; we then infer
from~\eqref{eq:exts} and~\eqref{eq:extmu} that we have $k=0$ and $B_1=B_2$
in \eqref{eq:extremal}, so that $M$ is K{\"a}hler--Einstein by Proposition
\ref{prop:newbiextremal}.

If, in addition, $W^-=0$, we obtain \eqref{eq:cond1} and $g$ is then a
selfdual K{\"a}hler--Einstein surface, i.e., a complex space form.

If $W ^-$ does not vanish identically, on the open set where $W^-\neq 0$,
the Hermitian pair $(g,I)$ is invariant under the action of the Killing
vector fields $K_i$'s, as $I$ is determined by the eigenform of $W^-$
corresponding to its simple eigenvalue, cf.~\cite{derdzinski,AG1}; then,
$\kappa$ (a constant multiple of the simple eigenvalue of $W^-$), the
square-norm $|\theta|^2$ of the Lee form $\theta$ of the pair $(g, I)$ as
well as $\delta \theta$ are also invariant under the action of $K _i$'s,
$i = 1, 2, 3$; we then have $d|\theta|^2 \wedge d\kappa =0$; by using
\eqref{eq:theta} and Lemma \ref{lem:curv}, we can check that this
implies that $F$ and $G$ satisfy \eqref{eq:cond2}; by Proposition
\ref{prop:newextremal}, the corresponding \biax/ K{\"a}hler surface is
Ricci-flat and $\kappa = - \frac{2 (C _1 - C _2)}{(\xi - \eta) ^3}$; it
then follows from \eqref{eq:theta} that $\theta = \frac{1}{3} \, d\ln
|\kappa|$; in particular, $d|\theta|^2 \wedge \theta =0$; by \cite[Theorem
1]{AG2}, this implies that the metric is of cohomogeneity one. Since
$\kappa$ is non-zero, it is not constant, and so the metric is not
homogeneous.

Evidently the two cases overlap when $C_1=C_2$ and $A=\ell=0$, in which case
$g$ is flat.
\end{proof}

\begin{rem} {\rm In fact we do not need to assume {\it a priori} that $M$ is
extremal in the above proof, as long as we assume that the additional
symmetry preserves $\varphi$, hence $\xi-\eta$. Then $ds$ and $d\xi-d\eta$
are collinear, from which it is easy to deduce that $M$ is extremal.}
\end{rem}

Let us now collect the results we have established so far about weakly
selfdual K{\"a}hler surfaces.

\begin{theo} Let $(M, g, J,\omega)$ be a weakly selfdual K{\"a}hler surface.
Denote by $s$, $\lambda$, and $p = \bigl(\frac{s}{2} + \lambda\bigr)
\bigl(\frac{s}{2} - \lambda\bigr)$, the \textup(normalized\textup) scalar
curvature, the positive eigenvalue of the trace-free Ricci tensor $\Ric
_0$, and the pfaffian of the normalized Ricci tensor, respectively. Then:

\smallskip\noindent{\rm (i)}
$K _1 := J {\rm grad} _g s$ and $K _2 := J {\rm grad} _g
p$ are commuting Killing vector fields, and on any simply connected
open subset where $K_1$ and $K_2$ are linearly independent, the
functions $\xi := \frac{s}{2} + \lambda$, $\eta :=
\frac{s}{2} - \lambda$, $t$, $z$ form a globally defined coordinate system
with respect to which the K{\"a}hler structure $(g,J,\omega)$ is
\begin{gather} \label{eq:newgbis} \begin{split}
g &= (\xi - \eta) \, \Bigl ( \frac{d \xi ^2}{F (\xi)}
- \frac{d \eta ^2}{G (\eta)} \Bigr)\\
&\qquad + \frac{1}{\xi - \eta} \bigl ( F (\xi) (dt + \eta dz) ^2
- G (\eta) (dt + \xi dz) ^2 \bigr),
\end{split}\\
\label{eq:newJbis} \begin{split}
J d \xi &= \frac{F (\xi)}{\xi - \eta} (dt + \eta\, dz), \qquad
J dt = -\frac{\xi \,d \xi}{F (\xi)} - \frac{\eta \,d \eta}{G (\eta)}, \\
J d \eta &= \frac{G (\eta)}{\eta - \xi} (dt + \xi\, dz), \qquad
J dz = \frac{d \xi}{F (\xi)} + \frac{d \eta}{G (\eta)},
\end{split}\\
\label{eq:newomegabis}
\omega=d\xi\wedge(dt+\eta\, dz)+d\eta\wedge (dt +\xi\, dz),
\end{gather}
where
\begin{align} \label{eq:Fter}
F (x) &= k x ^4 + \ell x ^3 + A x ^2 + B x + C _1,\\
\label{eq:Gter}
G (x) &= k x ^4 + \ell x ^3 + A x ^2 + B x + C _2.
\end{align}

\smallskip\noindent{\rm (ii)}
Conversely, each almost K{\"a}hler structure $(g, J,\omega)$ described
by {\rm \eqref{eq:newgbis}--\eqref{eq:Gter}} is K{\"a}hler and weakly
selfdual with
\begin{equation*}
s = - 2 k (\xi + \eta) -  \ell,  \qquad
p = 4 k ^2 \xi \eta + k \ell (\xi + \eta) + \frac{\ell ^2}{4},
\end{equation*}
so that $K _1 = \frac{\partial}{\partial t}$ and
$K _2 = \frac{\partial}{\partial z}$.

\smallskip\noindent{\rm (iii)}
The K{\"a}hler structure described by {\rm \eqref{eq:newgbis}--\eqref{eq:Gter}}
is selfdual if and only if $C_1=C_2$.
\end{theo}

In the selfdual case, we recover the general expression found by
Bryant~\cite[Section 4.3.2]{bryant} depending on an arbitrary
polynomial of degree $4$.

\section{K{\"a}hler surfaces of Calabi type}\label{s:coh1}

\subsection{Hamiltonian $2$-forms and the Calabi construction}

In this section we classify weakly selfdual K{\"a}hler surfaces of nowhere
constant scalar curvature $s$, but for which $p$ and $s$ are not
independent.  As in the previous section (when we assumed $p$ and $s$ were
independent) we do this by finding an explicit formula for a K{\"a}hler
surface $(M,g,J,\omega)$ with a \biham/ $2$-form
$\varphi=\varphi_0+\frac32\sigma\,\omega$ such that $K_1:=J{\rm grad}_g
\sigma$ has no zero, but $K_2:= J{\rm grad}_g\pi=b K_1$, where $\pi$ is the
pfaffian of $\tilde\varphi$ and $b$ is (necessarily) constant.

The general theory of K{\"a}hler surfaces with a \biham/ $2$-form,
described in Section \ref{s:wsdk}, applies equally to this case. In
particular, since $K_1$ has no zero, we may still write $\pi=\xi\eta$ and
$\sigma=\xi+\eta$ for the trace and pfaffian of $\tilde\varphi$, and $d\xi$
and $d\eta$ are orthogonal by Proposition~\ref{prop:dxi-deta}.

On the other hand, the constructions of Section \ref{s:baks} definitely
fail, because $K_1$ and $K_2$ are no longer independent: $\pi$ is an affine
function of $\sigma$, so $\xi$ and $\eta$ are not independent
functions. Therefore, $d\xi$ and $d\eta$, in addition to being orthogonal,
are collinear! This is not a contradiction: we deduce that either $\xi$ or
$\eta$ is constant. Since $\pi=\xi(\sigma-\xi)=\eta(\sigma-\eta)$, this
constant is the constant $b$ above, so that $\pi = b (\sigma - b)$ and
$\lambda = \pm \frac12(\sigma-2b)$.

We observe that $K _1$ is an eigenvector of $-\varphi_0\circ J$, for the
eigenvalue $\lambda$, and the conformally K{\"a}hler anti-selfdual complex
structure $I$ is characterized as follows: $I$ coincides with $J$ on the
distribution generated by $K _1$ and $J K _1$, but with $- J$ on the
orthogonal distribution. Hence we are in the following situation.

\begin{defi} A K{\"a}hler surface $(M,g,J,\omega)$ is said to be
\emph{of Calabi type} if it admits a non-vanishing hamiltonian Killing vector
field $K$ such the almost-Hermitian pair $(g,I)$---with $I$ equal to $J$ on
the distribution spanned by $K$ and $JK$, but $-J$ on the orthogonal
distribution---is conformally K{\"a}hler.
\end{defi}

It is straightforward to obtain an explicit formula for K{\"a}hler metrics
of Calabi type, using the LeBrun form of a K{\"a}hler metric with a
hamiltonian Killing vector field~\cite{lebrun}.

\begin{prop}\label{prop:monax} Let $(M,g,J,\omega)$ be a K{\"a}hler surface
of Calabi type. Then the K{\"a}hler structure is given locally by
\begin{align}\label{eq:mono-g}
g&=(a z - b)g_\Sigma+ w(z) dz^2  + w(z)^{-1} (dt+\alpha)^2,\\
\omega&= (a z-b)\,\omega_\Sigma+dz\wedge (dt+\alpha),
\label{eq:mono-om}\end{align}
where $z$ is the momentum map of the Killing vector field $K$, $t$ is a
function on $M$ with $dt(K)=1$, $w$ is a function of one variable,
$g_\Sigma$ is a metric on $2$-manifold $\Sigma$ with area form
$\omega_\Sigma$, $\alpha$ is a $1$-form on $\Sigma$ with
$d\alpha=a\,\omega_\Sigma$, and $a,b$ are constant.  Conversely equations
\mbox{\eqref{eq:mono-g}--\eqref{eq:mono-om}} define a K{\"a}hler
structure of Calabi type with $K=\partial/\partial t$, for any
$g_\Sigma$ and $V$.
\end{prop}
\begin{proof}
The proof follows LeBrun's description \cite{lebrun} of K{\"a}hler metrics
with a hamiltonian Killing vector field $K$. Supposing first that
$(g,J,\omega)$ is only almost Hermitian, note that $K - i JK$ is a
holomorphic vector field, so that the complex quotient is locally a Riemann
surface $\Sigma$. Introducing a local holomorphic coordinate $x+iy$
on $\Sigma$, we may write
\begin{align*} \label{eq:lebrunform}
g&= e ^u w(dx^2 + dy^2) + w \,dz^2 + w^{-1} (d t  + \alpha)^2,\\
&\quad J dx = dy,\qquad J dz = w^{-1} (dt+\alpha),\\
\omega &= e^u w \,dx\wedge dy+dz\wedge (dt+\alpha).
\end{align*}
where $dt(K)=1$, $\alpha$ is an invariant $1$-form with $\alpha(K)=0$, and
$u,w$ are functions of $x,y,z$. The almost Hermitian structure $I$ is
given by
\begin{equation*}
I dx = -dy,\qquad I dz = w^{-1} (dt+\alpha),
\end{equation*}
with K{\"a}hler form
\begin{equation*}
\omega_I = -e^u w\, dx\wedge dy + dz\wedge (dt+\alpha).
\end{equation*}
We now impose the condition that
$(g,J,\omega)$ and $(\bar
g=\lambda^{-2}g,I,\bar\omega=\lambda^{-2}\omega_I)$ are K{\"a}hler for some
non-vanishing function $\lambda$. Now $d\omega=0$ if and only if
\begin{equation}\label{eq:domJ=0}
(e^u w)_z dz\wedge dx\wedge dy=dz\wedge d\alpha,
\end{equation}
while $\lambda^{-2}\omega_I$ is closed if and only if
\begin{multline}\label{eq:omIcK}
\lambda\, dz\wedge d\alpha
+\bigl((e^u w)_z\lambda-2\lambda_z e^u w\bigr) dz\wedge dx\wedge dy\\
+2\lambda_x dx\wedge dz\wedge (dt+\alpha)
+2\lambda_y dy\wedge dz\wedge (dt+\alpha)
+2\lambda_t dt\wedge \omega_I=0.
\end{multline}
In the presence of~\eqref{eq:domJ=0}, \eqref{eq:omIcK} is equivalent to
\begin{equation}\label{eq:omIcK2}\begin{split}
(e^u w)_z\lambda=\lambda_z e^u w,\\
\lambda_x=\lambda_y=\lambda_t=0,
\end{split}\end{equation}
which holds if and only if $\lambda=\lambda(z)$ and $e^uw =h\lambda$ for
some function $h(x,y)$.

Let $\theta$ be the complex $1$-form $w\,dz+i (dt+\alpha)$. Then, since
$dx\pm i dy$ is closed, the complex structures $I$ and $J$ are integrable if
and only if $d\theta$ belongs to the ideals generated by $\{\theta, dx-i
dy\}$ and by $\{\theta,dx+i dy\}$ respectively.  Since
$d\theta\bigl(\frac\partial{\partial t},\cdot\bigr)=0$, these conditions
force $(dx-i dy)\wedge d\theta$ and $(dx+i dy)\wedge d\theta$ to vanish.
Hence $I$ is integrable if and only if
\begin{equation}\label{eq:Iint}
d\alpha=-w_x dy\wedge dz +w_y dx\wedge dz+f dx\wedge dy
\end{equation}
while $J$ is integrable if and only if
\begin{equation}
d\alpha=w_x dy\wedge dz -w_y dx\wedge dz+f dx\wedge dy;
\end{equation}
here $f$ is an arbitrary function.  Hence $I$ and $J$ are both integrable if
and only if
\begin{equation}\label{eq:IJint}
w_x=w_y=0 \qquad\text{and}\qquad d\alpha=f dx\wedge dy
\end{equation}
and $f$ is necessarily a function of $x,y$ only.

Putting together~\eqref{eq:domJ=0}, \eqref{eq:omIcK2}, and
\eqref{eq:IJint}, we see that $(g,J,\omega)$ is of Calabi type, with Killing
vector field $K$, if and only if $e^u w=h(x,y)\lambda(z)$, with
$d\alpha=h(x,y)\lambda_z$, $\lambda_{zz}=0$ and $w_x=w_y=0$,
so that $\lambda=a z-b$ for constants $a,b$ and $w=w(z)$.

Using the freedom in the choice of $t$, we may then assume $\alpha$ is a
$1$-form on $\Sigma$, while $g_\Sigma= h(x,y)(dx^2+dy^2)$ is a metric on
$\Sigma$, and the result follows.
\end{proof}

This Proposition shows that K{\"a}hler metrics of Calabi type are
essentially the same as metrics arising from the well-known \emph{Calabi
construction} \cite{calabi1} of metrics on the total space of a Hermitian
line bundle over a Riemann surface. In this interpretation, the Killing
vector field $K$ generates the natural circle action on the line bundle and
the K{\"a}hler form is
\begin{equation}
\omega_\Sigma + d J df
\end{equation}
for a function $f$ of the fibre norm $r$. Since $-d J d\log r$ is the
curvature of the line bundle, which is basic, the momentum map $z$
of $K$ is also a function of $r$. Hence we may locally view $f$
as a function of $z$, so that, if we write $Jdz=dt+\alpha$ where
$dt(K)=1$ and $\alpha$ is basic,
\begin{equation}
d J df = \Bigl(\frac{f'(z)}{w(z)}\Bigr)_{\!z} dz\wedge (dt+\alpha)
+\frac{f'(z)}{w(z)}d\alpha.
\end{equation}
Therefore, in Proposition~\ref{prop:monax}, setting $b=-1$ without loss of
generality, we have $f'(z)=z w(z)$.

The metrics of Proposition~\ref{prop:monax} certainly admit a
\biham/ $2$-form, namely $\varphi=(az-b)\omega_I+3a z \omega$. Hence
$\sigma=2az$ and $\lambda=az-b$, so that $\xi=2az-b$ and $\eta=b$. The
hamiltonian Killing vector fields associated to $\varphi$ both vanish when
$a=0$. On the other hand, for $a$ nonzero, we can use the freedom in the
choice of $z$ to set $a=1$ and $b=0$.

\begin{theo} A K{\"a}hler surface is of Calabi type if and only if either:
\begin{enumerate}
\item it is locally a K{\"a}hler product of two Riemann surfaces, one of
which admits a Killing vector field; or

\item it admits a \biham/ $2$-form whose associated Killing vector
fields are dependent but not both zero.
\end{enumerate}
The K{\"a}hler structure is then given explicitly by the Calabi
construction~\eqref{eq:mono-g}--\eqref{eq:mono-om}: in case {\rm (i)}
$a=0$, while in case {\rm (ii)} we may take $a=1$, $b=0$ without loss of
generality.
\end{theo}

\subsection{Weakly selfdual K{\"a}hler surfaces of Calabi type}

We begin this section by computing the curvature of a K{\"a}hler
surface of Calabi type which is not a local K{\"a}hler product. Therefore
we set $a=1$, $b=0$, and write $w(z)=z/V(z)$, so that the K{\"a}hler
structure is
\begin{align}\label{eq:gen-g}
g &= z g_\Sigma + \frac z{V(z)} dz^2  +\frac{V(z)}z (dt+\alpha)^2,\\
\omega&= z\,\omega_\Sigma+dz\wedge (dt+\alpha),
\label{eq:gen-om}\end{align}

As with \biax/ K{\"a}hler surfaces, the curvature is entirely determined by
the scalar curvature $s$ of $g$, the conformal scalar curvature $\kappa$ of
the Hermitian pair $(g,I)$, and the trace-free part $\rho _0$ of the Ricci
form of $(g, J)$.

\begin{lemma} \label{lem:mono-curv} For a non-product K{\"a}hler
surface $(M, g, J,\omega)$ of Calabi type, given
by~\eqref{eq:gen-g}--\eqref{eq:gen-om},
$\rho_0$ is a multiple $\mu$ of the K{\"a}hler form $\omega_I$ of the
Hermitian pair $(g, I)$, and $\mu,s,\kappa$ are given by
\begin{align}
\label{eq:mono-mu}
\mu &= -\frac1{4z} \biggl(s_\Sigma+\Bigl(\frac{V_z}{z^2} \Bigr)_{\!z}
z^2\biggr),\\
\label{eq:mono-scal}
s &=\frac{s_\Sigma-V_{zz}}{6z},\\
\label{eq:mono-kap}
\kappa &= \frac1{6z}\biggl(s_\Sigma-z^2
\Bigl(z^2\Bigl(\frac{V}{z^4}\Bigr)_{\!z}\Bigr)_{\!z}\biggr),
\end{align}
where $s_\Sigma$ is the scalar curvature of $\Sigma$.  In particular, on
the open subset of $M$ where $\mu$ has no zero, the anti-selfdual
almost-complex structure determined by $\rho _0$ is equal to $I$.
\end{lemma}
\begin{proof} The Ricci form $\rho$ is given by
$\rho=\rho_\Sigma -\frac12 d J d \log V$. The first term is $\frac12
s_\Sigma\omega_\Sigma$, and we compute the second term as follows:
\begin{align*}
dJd\log V &= d\Bigl(\frac{V_z}{V}J dz\Bigr)=
d\Bigl(\frac{V_z}{z} (dt+\alpha)\Bigr)\\
&= \frac {V_{zz}}{z} dz\wedge(dt+\alpha)+V_z\Bigl(\frac1z\omega_\Sigma
-\frac1{z^2} dz\wedge(dt+\alpha)\Bigr).
\end{align*}
Evidently we may write this as a linear combination of $\omega$ and
$\omega_I$, and we readily obtain~\eqref{eq:mono-mu} and~\eqref{eq:mono-scal}.

The conformal scalar curvature is most easily computed by noticing that
the conformal K{\"a}hler metric $\bar g= z^{-2}g$ is also of Calabi type,
with $\bar z=1/z$ and $\bar V(\bar z)= \bar z^4 V(1/\bar z)=V(z)/z^4$.
Hence its scalar curvature is
\begin{equation*}
\frac{s_\Sigma-\bar V_{\bar z\bar z}}{6\bar z}
= \frac{z}{6}
\biggl( s_\Sigma-z^2\Bigl(z^2\Bigl(\frac{V}{z^4}\Bigr)_{\!z}\Bigr)_{\!z}\biggr)
\end{equation*}
from which~\eqref{eq:mono-kap} follows, since $\kappa=z^{-2}\bar s$.
\end{proof}

\begin{prop}\label{prop:monoextremal}
Let $M$ be a non-product K{\"a}hler surface $M$ of Calabi type,
with Killing vector field $K$. Then the scalar curvature of $M$ is
a momentum map for a multiple of $K$ if
and only if $g_\Sigma$ has constant curvature $k$ and $V$ is of the form
\begin{equation} \label{eq:newP}
V(z)= A_1 z^4 + A_2 z^3 + k z^2 + A_3 z + A_4.
\end{equation}
Any K{\"a}hler surface given by~\eqref{eq:gen-g}--\eqref{eq:gen-om}
is extremal, with Ricci form $\mu\,\omega_I+\frac32 s\,\omega$, where
\begin{align}\label{eq:mu-coh}
\mu&= - A_1 z + \frac{A_3}{2z^2}\\
s &=  - 2 A_1 z -  A_2.\label{eq:s-coh};\\
\intertext{also the conformal scalar curvature of $(g<I)$ is}
\kappa &= - \frac{A_3}{z^2}-\frac{2 A_4}{z^3}.
\label{eq:kappa-coh}
\end{align}
Hence:
\begin{enumerate}
\item $g$ has constant scalar curvature if and only if $A_1=0$;

\item $g$ is scalar-flat \textup(i.e., anti-selfdual\textup) if
and only if $A_1=A_2=0$.

\item $(g, J)$ is K{\"a}hler--Einstein if and only if
$A _1 = A _3 = 0$;

\item $g$ is weakly selfdual if and only if $A _3 = 0$.

\item $g$ is selfdual if and only if $A _3 = A _4 = 0$.

\item $(g, J)$ is bi-extremal if and only if $g$ is weakly
selfdual;

\item the Bach tensor of $g$ vanishes if and only if
$4 A _1 A _4 - A _2A _3 = 0$.
\end{enumerate}\end{prop}
\begin{proof} From~\eqref{eq:mono-scal}, we have
\begin{equation*}
6zs +V_zz =s_\Sigma.
\end{equation*}
If $s$ is an affine function of $z$, then both sides of this equation must
be constant. Hence $s_\Sigma=2k$ and $V$ must be a quartic in $z$ with
quadratic term $k z^2$. The formulae for $\mu, s$ and $\kappa$ are immediate,
as are (i)--(v). (For (iv) we use the fact that $I$ is integrable
and $z^{-2}\omega_I$ is closed, so that weak selfduality is equivalent
to the equation $ds=2 d\mu$.)

(vi)  The pfaffian $p$ of the normalized Ricci form is given by
\begin{equation}
p = \Bigl(- 2 A_1 z - \frac12 A_2+\frac{A_3}{2z^2}\Bigr)\Bigl(-\frac12 A_2
-\frac{A_3}{2z^2}\Bigr)
\end{equation}
in particular, $p$ is a rational function of $z$; since $z$ is a
holomorphic potential, $p$ is a holomorphic potential if and only if it is
an affine function of $z$, if and only if $A _3 = 0$; $p$ is then equal
to $- \frac12 A _2 (s + \frac12 A _2)$.

(vii) For any extremal K{\"a}hler surface the Bach tensor is $J$-invariant
and the corresponding anti-selfdual $2$-form, is expressed by
\eqref{eq:Bachform} which yields
\begin{equation*}
\tilde{B} = \frac{(4 A _1 A _4 - A _2 A _3)}{z ^2} \, \omega _I.
\end{equation*}
\end{proof}

This family of extremal K{\"a}hler metrics has been considered in many
places.  In particular, it includes the extremal K{\"a}hler metrics of
cohomogeneity one under $U (2)$ constructed by Calabi in \cite{calabi1};
more generally, it turns out that these metrics all have cohomogeneity one
under a (local) action of a four-dimensional Lie group, locally isomorphic
to a central extension of the isometry group of a surface of constant
curvature $k$.  We refer to these metrics as {\it extremal K{\"a}hler
surfaces of Calabi type} and briefly recall how they may be realized as
diagonal Bianchi metrics, of class IX, VIII or II, according to whether
$k$ is positive, negative or zero.

Up to rescaling, we can---and will---assume that $k = \varepsilon$, with
$\varepsilon = 1, 0$ or $ -1$. As is well known, we can now write
$dt+\alpha=\sigma_3$ and introduce $t$-dependent $1$-forms
$\sigma_1,\sigma_2$ on $\Sigma$ with
$g _{\Sigma} = \sigma _1^{\,2} + \sigma _2^{\,2}$,
$\sigma_1\wedge\sigma_2=\omega_\Sigma$ and
\begin{equation} \label{eq:sigmas}
d \sigma _1 = \varepsilon \sigma _2 \wedge \sigma _3, \ \
d \sigma _2 = \varepsilon \sigma _3 \wedge \sigma _1, \ \
d \sigma _3 = \sigma _1 \wedge \sigma _2.
\end{equation}

By substituting $\sigma _1, \sigma _2, \sigma _3$ in \eqref{eq:gen-g}
\begin{equation} \label{eq:bianchi}
g = \frac{z}{V(z)} dz^2 + z ( \sigma _1^{\,2} + \sigma _2^{\,2} )
+ \frac{V(z)}{z} \, \sigma _3^{\,2};
\end{equation}
the complex structure $J$ is determined by
\begin{equation} \label{eq:Jbianchi}
J \sigma _1 = \sigma _2, \qquad
J d z = \frac{V (z)}{z} \, \sigma _3
\end{equation}
while the K{\"a}hler form $\omega$ is
\begin{equation} \label{eq:Jomega}
\omega = d z \wedge \sigma _3 + z\sigma _1 \wedge \sigma _2 = d (z \sigma _3).
\end{equation}

We here recognize {\it bi-axial diagonal Bianchi metrics}
of class IX, VIII or II, according as $\varepsilon$ is equal $1$, $-1$ or
$0$; these admit a cohomogeneity one local action of $SU (2)$ if $\varepsilon
= 1$, of $SU(1,1)$ if $\varepsilon = - 1$, of the Heisenberg group Nil if
$\varepsilon = 0$, and the orbits are level sets of $z$. This can be seen as
follows: denote by $(Z _1, Z _2 , Z _3 = K _1)$ the triplet of vector fields
determined by $\sigma _i (Z _j) = \delta _{ij}$ and $d z (Z _i) = 0$, $i ,
j = 1, 2, 3$, where $\delta _{ij}$ is the Kronecker symbol.  Then, for each
value of $z$, $Z _1, Z _2, Z _3$ are tangent to the corresponding orbit $M
_{z}$ and generate a Lie algebra isomorphic to $\mathfrak{su} (2)$,
$\mathfrak{su} (1, 1)$ or $\mathfrak{nil} ^3$ according as $\varepsilon$ is
equal to $1$, $-1$ or $0$; therefore, each orbit can be locally identified
to the corresponding Lie group and we can locally construct a new triple of
independent vector fields, $(\tilde{Z} _1, \tilde{Z} _2, \tilde{Z} _3)$, such
that $[\tilde{Z} _i, Z _j] = 0$ and $d z (\tilde{Z} _i) = 0$, $i,j = 1,2,
3$ (for each orbit, if $(Z _1,Z _2, Z _3)$ is a basis of
left-invariant vector fields, $(\tilde{Z} _1, \tilde{Z} _2, \tilde{Z} _3)$
is a basis of right-invariant vector fields); $\tilde{Z} _1,
\tilde{Z} _2, \tilde{Z} _3$ are clearly Killing with respect to $g$ and all
commute with $K _1$.

If $\varepsilon \neq 0$, $K _1, \tilde{Z} _1, \tilde{Z} _2, \tilde{Z} _3$
generate a $4$-dimensional Lie algebra, corresponding to a local action of
$U (2)$ if $\varepsilon = 1$, of $U (1, 1)$ if $\varepsilon = -1$.  If
$\varepsilon = 0$, $\tilde{Z} _3$ equals $K _1$, up to a constant factor;
on the other hand, we get an additional Killing vector field $\tilde{K} _1$
generated by the rotations around the origin in the Euclidean $2$-plane
$\E^2$ of $x, y$; then, $\tilde{K} _1, K _1, \tilde{Z _1}, \tilde{Z _2}$
generate a $4$-dimensional Lie algebra, say ${\mathfrak g}$, corresponding
to a local action of the group, $G$, obtained by forming the semi-direct
product of ${\rm Nil}$ by $S ^1$ for the natural action of $S ^1$ on ${\rm
Nil}$ by (outer) automorphisms; the centre of $G = {\rm Nil} \ltimes S ^1$
coincides with the centre of ${\rm Nil}$; the latter is one dimensional
again and the quotient of $G$ by its center is isomorphic to ${\rm Isom}
({\E^2})$; in other words, $G$ is isomorphic to a (non-trivial)
one-dimensional central extension of ${\rm Isom} ({\E^2})$.

If we concentrate our attention on weakly selfdual K{\"a}hler surfaces, we
readily infer from the foregoing:

\begin{theo}\label{th:monowsd}
Let $(M,g, J)$ be a weakly selfdual K{\"a}hler surfaces.  Denote by $s$
the scalar curvature and by $p$ the pfaffian of the normalized Ricci form;
let $K _1 = J {\rm grad} _g s$ and $K _2 = J {\rm grad} _g p$ be the
associated Killing vector fields and assume that $K _1$ has no zero and
that $K _2 = b K _1$, where $b$ is a real constant \textup(possibly
zero\textup).

\smallskip\noindent{\rm (i)}
Then $(M, g, J)$ admits a local action of cohomogeneity one of $G = U (2)$,
$U (1, 1)$ or ${\rm Nil} \ltimes S ^1$ and is locally isomorphic to a
diagonal Bianchi metric of class IX, VIII or II respectively.

More precisely, let $\sigma _1, \sigma _2, \sigma _3$ denote the
\textup(local\textup) $1$-forms on $M$ induced by this action,
corresponding to a triple of $G$-invariant $1$-forms on $G$, so that
$d \sigma _3 = \sigma _1 \wedge \sigma _2$,
$d \sigma _2 = \varepsilon \sigma _3 \wedge \sigma _1$,
$d \sigma _1 = \varepsilon \sigma _2 \wedge \sigma _3$,
where $\varepsilon = 1, -1$ or $0$ according as $G = U (2)$, $U (1, 1)$ or
${\rm Nil} \ltimes S ^1$; then, the K{\"a}hler structure $(g, J)$
can be put in the form {\rm \eqref{eq:bianchi}--\eqref{eq:Jbianchi}}, where
$z$ is an affine function of $s$ and $V (z)$ is of the form
\begin{equation} \label{eq:newwaP}
V (z) = A _1 z ^4 + A _2 z ^3 + \varepsilon z ^2 + A_4.
\end{equation}

\smallskip\noindent{\rm (ii)} Conversely, each K{\"a}hler surface of the form
{\rm \eqref{eq:bianchi}--\eqref{eq:Jbianchi}}, with $V$ given by
\eqref{eq:newwaP} is weakly selfdual.

\smallskip\noindent{\rm (iii)} This K{\"a}hler structure is selfdual if and
only if $A_4 = 0$.
\end{theo}

\begin{rem} \label{eq:remA4C} {\rm By substituting $A _3 = 0$ in
\eqref{eq:mu-coh} and \eqref{eq:kappa-coh}, we readily infer that that $A
_4$ and the constant $c$ appearing in \eqref{eq:C} and \eqref{eq:Ckappa}
are linked together by
\begin{equation} c = 2 A _1 ^3 A _4; \end{equation}
moreover, $A _1 \neq 0$, as $s$ is non-constant, so that $A _4 = 0$ if and
only if $c = 0$; as we already know, both conditions are equivalent to $g$
being selfdual.}
\end{rem}

\subsection{The weakly selfdual K{\"a}hler metric on the first
Hirzebruch surface}

We close this section by providing an example of a {\it compact} weakly
selfdual K{\"a}hler surface $(M, g, J)$; this belongs to the family of
extremal metrics constructed in \cite{calabi1} by E. Calabi on the first
{\it Hirzebruch surface} $F _1$, viewed as the compactification of the
total space of the tautological line bundle $L = {\mathcal O} (-1)$ over
the complex projective line ${\C} P ^1 = {\mathbb P} ({\C} ^2)$ obtained by
adding a {\it section at infinity}, say $C _{\infty}$; then, the zero
section, $C _0$, has self-intersection $-1$ and $F _1$ can also be
considered as a blown up of the projective plane ${\mathbb C} P ^2$ at some
point, with exceptional divisor $C _0$, and $F _1 - C _0 = {\C} ^2 - \{ (0,
0) \}$.

We here recall the main features of Calabi's construction from
\cite{calabi1}.  Let $r$ be the usual norm in ${\C} ^2 $ and, for
convenience, introduce the function $t$ defined on ${\C} ^2 - \{ (0, 0) \}$
by $e ^t = r ^2$; the group $U (2)$ naturally acts on $F _1$, by preserving
$C _0$ and $C _{\infty}$; by \cite{calabi2}, the connected group of
isometries, ${\rm Isom}_0(M,g)$, of any compact extremal K{\"a}hler surface
$(M, g, J)$ is a maximal compact subgroup of the connected group of
holomorphic transformations ${\rm Aut}_0(M,J)$; it follows that any
extremal metric on $F _1$ is isometric to a $U (2)$-invariant metric, in
particular has a globally defined potential on the open set
$M _0 = F _1 - C _0 = {\C} ^2 - \{ (0, 0) \}$ of the form $u = u (t)$; the
K{\"a}hler form is thus given by $\omega = \psi (t) d d ^c t + \psi ' (t)
dt \wedge d ^c t$, where $\psi (t) := u' (t)$; conversely, any such
K{\"a}hler metric extends to $F _1$ if and only if $\psi (t)$ extends to a
$C ^{\infty}$ function of $e ^ t$ in the neighbourhood of $t = - \infty$
and to a $C ^{\infty}$ function of $e ^ {- t}$ in the neighbourhood of $t =
\infty$; in particular, $\psi$ has a limit $a$ when $t \to - \infty$ and a
limit, $b$, when $t \to \infty$,with $0 < a < b$; the corresponding
K{\"a}hler class $[\omega]$ is then equal to $4 \pi (- a [C _0] + b [C
_{\infty}])$, where $[C _0]$ and $[C _{\infty}]$ denotes the
Poincar{\'e}-dual of $C _0$ and similarly for $[C _{\infty}]$; conversely,
each K{\"a}hler class on $F _1$ is of this form,for some pair $0 < a < b$;
it is easily checked that the Ricci for is given by $\rho = d d ^c v$, with
$v = t - \frac{1}{2} \log{\psi} - \frac{1}{2} \log{\psi '}$, so that $s = 2
\big ( \frac{v '}{\psi} + \frac{v ''}{\psi '} \big )$; on the other hand,
such a metric is extremal if and only if $s$ is an affine function of
$\psi$; by easy successive partial integrations, we infer that this metric
is extremal if and only if $\psi$ satisfies the following differential
relation: $\psi \psi ' = V (\psi)$, where $V$ is a polynomial of the form
\eqref{eq:newP}; moreover, the extremal metric is actually defined on $F _1$
and has its K{\"a}hler class parameterized by the pair $(a, b)$ as above if
and only if the coefficients of $V$ are given by
\begin{equation} \label{eq:calabiP} \begin{split}
& A _1 = \frac{- 2a}{(b - a) (a ^2 + 4 ab + b ^2)}, \\
& A _2 = \frac{(3 a ^2 - b ^2)}{(b - a) (a ^2 + 4 ab + b ^2)}, \\
& A _3 = \frac{a b (3 a ^2 - b ^2)}{(b - a) (a ^2 + 4 ab + b ^2)}, \\
& A _4 = \frac{- 2a ^3 b ^2}{(b - a) (a ^2 + 4 ab + b ^2)};
\end{split} \end{equation}
we thus get an extremal K{\"a}hler metric, say $g _{(a, b)}$, in each
K{\"a}hler class of $F _1$ (we have a similar construction for the other
Hirzebruch surfaces $F _k$ \cite{calabi1}; the coefficients of $V$ are then
given by \eqref{eq:parak} below).

All these metrics are of cohomogeneity one under the action of $U (2)$ and
can also be put in the form \eqref{eq:bianchi}--\eqref{eq:Jbianchi}, with
the same parameter $\psi$, the same polynomial $V (\psi)$, $d ^c t = \sigma
_3$ and $d d ^c t = \sigma _1 \wedge \sigma _2$; in particular, according
to Proposition \ref{prop:monoextremal}, $g _{(a, b)}$ if weakly selfdual
if and only if $A _3 = 0$, and this happens if and only if $a$ and $b$ are
related by
\begin{equation} \label{eq:chernwa} b ^2 = 3 a ^2; \end{equation}
in this case, $A _2 = 0$ as well, meaning that $p = 0$, and
$V(\psi) = A _1 \psi ^4 + \psi ^2 + A _4$ is a function of $\psi ^2$;
in particular, $\psi \psi ' = V (\psi)$ is easily integrated into
\begin{equation} \psi =
a \Bigl( \frac{1 + 3 e ^{t + t_0}}{1 + e ^{t + t_0}} \Bigr) ^{\frac{1}{2}} =
a \Bigl( \frac{1 + 3 e ^{t _0}  r ^2}{1 + e ^{t_0}r ^2}\Bigr) ^{\frac{1}{2}},
\end{equation}
where $t _0$ is a constant; the latter can be made equal to zero without
loss of generality by a mere translation of the parameter $t$, i.e. a
rescaling of $r$; then, up to rescaling, the K{\"a}hler form is given by
\begin{equation} \label{eq:calabiwa} \omega =
\frac{(1 + 3 e ^{t}) ^{\frac{1}{2}}}{(1 + e ^{t}) ^{\frac{1}{2}}} d d ^c t
+ \frac{e ^t}{(1 + e ^t) ^{\frac{3}{2}} (1 + 3 e ^t) ^{\frac{1}{2}}} dt
\wedge d ^c t.
\end{equation}

In the sequel, the K{\"a}hler metric given by \eqref{eq:calabiwa} will be
referred to as the {\it Calabi weakly selfdual K{\"a}hler metric} of
$F _1$. It may be noticed that, according to \eqref{eq:calabiP}, the scalar
curvature $s = - 2 A _1 \psi$ is non-constant and (strictly) positive.

In the next section we show that, conversely, each compact weakly selfdual
K{\"a}hler surface with non-constant scalar curvature is, up to rescaling,
isomorphic to the first Hirzebruch surface equipped with the Calabi weakly
selfdual K{\"a}hler metric.

\section{Compact weakly selfdual K{\"a}hler surfaces}\label{s:cwsd}

Compact {\it selfdual} K{\"a}hler surfaces have been described by
B.-Y. Chen in \cite{BYC}: these are locally symmetric, hence of constant
holomorphic sectional curvature or, locally the product of Riemann surfaces
of opposite constant curvature (see \cite{bryant} for a higher dimensional
generalization).

In \cite{jelonek}, W. Jelonek proved that compact real analytic weakly
selfdual K{\"a}hler surfaces are either K{\"a}hler--Einstein, or locally the
product of two Riemann surfaces of constant Gauss curvatures, or
biholomorphic to a ruled surface.

\smallskip
We show that the hypothesis of real analyticity can actually be removed and
that, except in the case when the scalar curvature is constant, the only
weakly selfdual K{\"a}hler ruled surface is the first Hirzebruch surface
$F _1$ equipped with the Calabi weakly selfdual K{\"a}hler metric, as
described in the preceding section (up to rescaling). More precisely
we have the following result.

\begin{theo} \label{theo3} Let $(M, g, J)$ be a compact weakly selfdual
K{\"a}hler surface. Then $(M,g,J)$ is either
\begin{enumerate}
\item  K{\"a}hler--Einstein; or
\item locally isomorphic to the product of two Riemann surfaces of
constant Gauss curvatures, or
\item up to rescaling, isomorphic to the first Hirzebruch surface $F_1$
equipped with a Calabi weakly selfdual K{\"a}hler metric.
\end{enumerate}
\end{theo}

\begin{proof}
By Proposition \ref{prop:tri}, we know that $g$ is either selfdual, hence,
by the above mentioned result of B.-Y. Chen, described by (i) or (ii), or
of constant scalar curvature, hence, again, described by (i) or (ii), or of
non constant scalar curvature.  In the latter case, the (negative)
K{\"a}hler structure $(\bar{g}, I)$ is globally defined; it then follows
from a result of Kotschick \cite{kotschick} that the signature of $M$ is
zero \cite{kotschick}; moreover, since the (real) holomorphic field $K_1 =
J {\rm grad}_g s$ has non-empty zero set, we know by \cite{kosniowski} that
the Kodaira dimension of the K{\"a}hler surface $(M,J)$ is $-\infty$,
hence $(M,J)$ is a ruled surface which is the projectivization
${\mathbb P}(E)$ of a rank 2 holomorphic vector bundle $E$ over a compact
complex curve
$\Sigma$
\cite{BPV}.

\smallskip
If $\Sigma = {\C P}^1$, $(M,J)$ is a {\it Hirzebruch surface} $F_k=
{\mathbb P}({\mathcal O} \oplus {\mathcal O}(-k))$, where $k$ is a positive
integer, or the product ${\mathbb CP}^1\times {\C P}^1$; the only extremal
K{\"a}hler metrics of ${\mathbb CP}^1\times {\C P}^1$ are the (symmetric)
product metrics, which are of constant scalar curvature; on the other hand,
any maximal compact subgroup of ${\rm Aut}_0(M,J)$ is conjugate to $U(2)$,
and therefore any extremal K{\"a}hler metric must be a cohomogeneity-one
$U(2)$ metric \cite{calabi2}, hence, locally of the form
\eqref{eq:bianchi} with
$\varepsilon = 1$ (cf.~the end of the preceding section for the case when
$k = 1$); as shown by E. Calabi in \cite{calabi1}, for each $k>0$, any
K{\"a}hler class of $F_k$ carries a unique extremal K{\"a}hler metric (up
to a re-parameterization); each one can be put in the form
\eqref{eq:bianchi}, where the polynomial $V$, in the notation of
\eqref{eq:newP}, is determined by
\begin{equation} \label{eq:parak} \begin{split}
& A _1 = \frac{(k + 1) a + (k - 1) b}{(b - a)(a ^2 + b ^2 + 4ab)}\\
& A _2 = \frac{(2 - k)b ^2 - (k + 2)a ^2}{(b - a)(a ^2 + b ^2 + 4ab)},\\
& A _3 = \frac{ab((2 - k)b ^2 - (k + 2)a ^2)}{(b - a)(a ^2 + b ^2 + 4ab)},\\
& A _4 = \frac{a ^2b ^2 ((k + 1)a + (k - 1)b)}{(b - a)(a ^2 + b ^2 + 4ab)},
\end{split} \end{equation}
where $0 < a < b$ are the parameters of the K{\"a}hler class; according to
Proposition \ref{prop:monoextremal}, $g$ is weakly selfdual precisely when
$A _3=0$; in the present situation, this is equivalent to $A _2=0$ and
happens if and only if $k=1$ and $\mu := \frac{a}{b} = \frac{1}{\sqrt{3}}$,
i.e. if $(M,g, J)$ is the first Hirzebruch surface equipped with a Calabi
weakly selfdual K{\"a}hler metric.

\smallskip
We now show that a compact ruled surface $(M, J)$ whose base $\Sigma$ is a
compact complex curve of genus ${\bf g}(\Sigma)$ at least 1 does not carry
weakly selfdual K{\"a}hler metrics of non-constant scalar curvature. We
thus assume that $(M, J) = {\mathbb P}(E)$ carries a weakly selfdual
K{\"a}hler metric of non-constant scalar curvature to get a contradiction.

Using an argument from
\cite{lebrun-singer,christina},  we first observe that the rank two vector
bundle $E$ splits as
$E ={\mathcal O} \oplus L$, where ${\mathcal O}$ stands for the trivial
holomorphic line bundle and $L$ is a holomorphic line bundle $L$ of degree
${\rm deg} (L) > 0$. Indeed, recall the already mentioned result of E.
Calabi \cite{calabi2} that the connected component of the isometry group
${\rm Isom}_0(M,g)$ is a maximal compact subgroup in ${\rm Aut}_0(M,J)$;
according to M. Maruyama
\cite{maruyama}, the group of automorphisms of ruled surfaces can be
described as follows: If ${\bf g}(\Sigma) \ge 1$, there exists an exact
sequence
\begin{equation}\label{eq:seq}
\{ {\bf 1} \} \to {\rm Aut}_{\Sigma}({\mathbb P} (E))
\to {\rm Aut}({\mathbb P} (E))\to {\rm Aut}(\Sigma),
\end{equation}
where ${\rm Aut}_{\Sigma}({\mathbb P} (E))$ denotes the group of relative
automorphisms of the bundle ${\mathbb P}(E)\to \Sigma$, and
${\rm Aut}(\Sigma)$ is the group of automorphisms of $\Sigma$ (of course,
${\rm Aut}(\Sigma)$ is finite if ${\bf g}(\Sigma) \ge 2$); on the other
hand, the (non-trivial) homomorphic vector field $\Xi _1 = K_1 - iJK_1$
whose real part is the Killing vector field $K_1= J {\rm grad}_g s$ has a
non-empty zero set; since $\Xi _1$ preserves the (unique) ruling
$M={\mathbb P}(E) \to \Sigma$, it projects onto a holomorphic vector field
on the base $\Sigma$; since $\Xi_1$ has at least one zero, the induced
vector field on $\Sigma$ vanishes; it follows that $\Xi _1$ is tangent to
the ${\C P}^1$ fibers (equivalently, $\Xi _1$ belongs to the Lie algebra of
${\rm Aut}_{\Sigma}({\mathbb P} (E))$); this shows that the kernel of the
group homomorphism $f: {\rm Isom}_0(M,g) \to {\rm Aut}_0(\Sigma)$ induced
by the exact sequence \eqref{eq:seq} is a non-trivial compact subgroup of
${\rm Aut}_{\Sigma}({\mathbb P} (E))$; one can therefore find an $S^1$ in
the connected component of the identity of ${\rm Aut}_{\Sigma}({\mathbb P}
(E))$; denote by $\Xi _0$ the induced holomorphic vector field, such that
the imaginary part $K_0= Im(\Xi_0)$ generates the $S^1$-action, whereas
$\Xi _0$ itself generates a ${\mathbb C}^{*}$-action; as a matter of fact,
$\Xi_0$ can be identified to a traceless holomorphic section of ${\rm
End}(E)$, say
${\bf s}$; note that ${\bf s}$ is of constant determinant; since
$K_0=Im(\Xi_0)$ generates a periodic $S^1$-action, ${\bf s}$ must be
diagonalizable; this shows that we have a holomorphic splitting of $E$
into eigensubbundles of
${\bf s}$; by twisting by a line bundle, we obtain the splitting $E
={\mathcal O} \oplus L$, where ${\rm deg} (L) \ge 0$; then, $\Xi_0$ is
nothing but the Euler vector field of $L$.  If the degree of $L$ is zero,
then any K{\"a}hler class contains a locally symmetric K{\"a}hler metric
\cite{wall}, so that any extremal K{\"a}hler metric on $(M,J)$ is of
constant scalar curvature \cite{calabi1}, a contradiction; we thus obtain a
splitting $E= {\mathcal O}\oplus L$ where $L$ is a holomorphic line bundle
of ${\rm deg} (L) > 0$.

As ${\rm Isom}_0(M,g)$ is a maximal compact subgroup in ${\rm Aut}_0(M,J)$,
we may assume \cite{christina} that (up to a biholomorphism) the metric $g$
is invariant under the fixed $S^1$ action generated by $K_0=Im(\Xi_0)$. For
any non-trivial Killing vector field, $K$, which arises from a real
holomorphic potential, the argument already used above shows that $\Xi=
K-iJK$ must be tangent to the fibers, and therefore $\Xi \wedge \Xi_0 =0$.
In particular, we get that $ K_0 = f K_1 + h JK_1, $ where $f,h$ are smooth
functions defined on an open dense subset of $M$ where $K_1 = J {\rm grad}_g
s \neq 0$. But $\g{K_0, JK_1} = - ds(K_0) = -{\mathcal L}_{K_0} s =0$, i.e.
$h=0$, and therefore $f$ is a constant. By rescaling the metric if necessary
we may assume therefore $K_1 = K_0= Im(\Xi_0)$.  Similarly, $K_2= J{\rm
  grad}_g p$ must be a constant multiple of $K_1$ and by Theorem
\ref{th:monowsd}, $g$ must be locally of cohomogeneity one, i.e., $g$ can be
written of the form \eqref{eq:bianchi} on an open dense subset of $M$.

Note that $M$ contains exactly two curves
fixed by the ${\mathbb C}^*$-action generated by $\Xi_0$, corresponding to
the zero and infinity sections, $C_0$ and $C_{\infty}$, of
$M={\mathbb P}({\mathcal O}\oplus L)$; moreover, the function $z$
appearing in \eqref{eq:bianchi} makes sense on the whole of $M$ as being a
momentum map of the corresponding $S^1$-action (up to multiplication by a
non-zero constant); it then follows that $z : M \to {\R}$ maps $M$ onto
an interval $[a,b]$, such that $z$ is regular on
$M - (z^{-1}(a)\cup z^{-1}(b))$; therefore, for any $t_0 \in (a,b)$,
$\Sigma =z^{-1}(t_0)/ S^1$, whereas $C_0 = z^{-1}(a)$ and
$C_{\infty}= z^{-1}(b)$ (see \cite{lebrun-singer}). By using an argument
from \cite[p.42]{lebrun-singer}, it is shown that $q = |K_1|^2$ is a smooth
function on $\Sigma \times [a,b]$, which satisfies the boundary conditions
\begin{equation}\label{eq:boundary}
q(.,a)=q(.,b)=0; \ (\frac{\partial}{\partial z} q)(.,a)
= - (\frac{\partial}{\partial z} q)(.,b) = k,
\end{equation}
where $k$ is a real constant; for Calabi's metrics \eqref{eq:bianchi} one
has $q = \frac{V(z)}{z}$, so that the equations \eqref{eq:boundary} read
$$V(a) = V(b)= 0; V'(a) = k a; V'(b) = -k b;$$
we thus obtain the following values for the coefficients $A _1, A _2,
A _3, A _4$ of $V$ (notations of \eqref{eq:newP}):
\begin{equation}\label{eq:parak1} \begin{split}
& A _1 =
\frac{ k (a + b) + \varepsilon (a - b)}{(b - a) (a ^2 + 4 a b + b ^2)},\\
& A _2 = \frac{ - k (a ^2  + b ^2) + 2 \varepsilon (b ^2 - a ^2)}
{(b - a) (a ^2 + 4 a b + b ^2)}, \\
& A _3 = \frac{ a b ( - k (a ^2 + b ^2) + 2 \varepsilon (b ^2 - a ^2))}
{(b - a) (a ^2 + 4 a b + b ^2)}, \\
& A _4 = \frac{a ^2 b ^2 ( k (a + b) + \varepsilon (a - b))}
{(b - a) (a ^2 + 4 a b + b ^2)};
\end{split} \end{equation}
according to Proposition \ref{prop:monoextremal} we have $A _3 =0$, and from
\eqref{eq:parak1} we also get $A _2 =0$; by \eqref{eq:s-coh} and
\eqref{eq:mu-coh} it follows that the Ricci tensor of $g$ has two distinct
eigenvalues, equal to $\frac{s}{6}$ and $\frac{s}{3}$ respectively; by
Proposition \ref{prop:tri}, $\Ric_0$ nowhere vanishes on $M$, meaning that
the scalar curvature $s$ nowhere vanishes as well; since $K_1$ is a
non-trivial Killing vector field, $s$ must be everywhere positive; this
shows that the first Chern class $c_1(M)$ of $(M,J)$ is positive, and
therefore $c_1^2(M) > 0$, a contradiction \cite{BPV}.
\end{proof}

\begin{rem} {\rm The case (ii) of Theorem \ref{theo3} includes in
particular ruled surfaces ${\mathbb P} (E)$ over a Riemann surface $\Sigma$
of genus $g \ge 1$ when the holomorphic vector bundle $E$ is stable or
polystable, cf.~e.g. \cite{wall}.
}
\end{rem}

\section{Weakly selfdual almost K{\"a}hler manifolds}

\subsection{The Matsumoto-Tanno identity for almost K{\"a}hler $4$-manifolds.}
Recall that an {\it almost K{\"a}hler} manifold is an almost Hermitian
manifold $(M,g,J,\omega)$ for which the K{\"a}hler $2$-form $\omega$ is
closed. The almost complex structure $J$ of an almost K{\"a}hler manifold is
not integrable in general; if it is, we obtain a K{\"a}hler manifold.

We would like to identify the Ricci tensor $\Ric$ of an almost K{\"a}hler
manifold with a $2$-form $\rho$, the Ricci form, as in the K{\"a}hler
case. However, only the $J$-invariant part of $\Ric$ defines a $2$-form,
whereas on a (non-integrable) almost K{\"a}hler manifold, the Ricci tensor is
not in general $J$-invariant. We shall therefore impose $J$-invariance as an
extra requirement.

Throughout this section we will always assume that $(M,g,J)$ is an almost
K{\"a}hler $4$-manifold whose Ricci tensor is $J$-invariant, i.e.,
$\Ric(J\cdot ,J\cdot) = \Ric(\cdot ,\cdot)$.  We then adopt the notations of
Section~\ref{s:wsdk}, and, in analogy with the K{\"a}hler case, we consider
the type $(1,1)$ {\it Ricci form}, $\rho$, of $(M,g,J)$ defined by
$\rho(\cdot, \cdot) = \Ric(J \cdot, \cdot)$; the anti-selfdual part of
$\rho$ is denoted by $\rho _0$.  It is a remarkable fact \cite{Dr} that even
though $J$ is not integrable, $\rho$ is still a closed (1,1)-form (although
it is no longer a representative of $\frac{1}{2\pi}c_1^{\mathbb R}$).  Using
this observation, the proof of Lemma \ref{lemma1} easily extends to the case
of almost K{\"a}hler $4$-manifolds with $J$-invariant Ricci tensor.

\begin{lemma} \label{ak-lemma1} For any almost K{\"a}hler $4$-manifold
with $J$-invariant Ricci tensor the identity {\rm (\ref{eq:bianchirho})} is
satisfied. In particular, the anti-selfdual Weyl tensor $W^-$ of such a
manifold is harmonic if and only if the Matsumoto-Tanno identity {\rm
(\ref{eq:MT})} is satisfied.
\end{lemma}
\begin{proof}
The proof follows the one given in the K{\"a}hler case, with slight
modifications in places where the non-integrability of $J$ must be taken
into account: the Cotton-York tensor is now written as
\begin{equation}\label{basicpart} \begin{split}
 C _{X, Y}(Z) =
& -\frac{1}{2}\Big((\nabla _X \rho)(Y, JZ) - (\nabla_Y \rho)(X, JZ) \Big) \\
& -\frac{1}{2}\Big(\rho(Y,(\nabla_X J)(Z)) - \rho(X,(\nabla_Y J)(Z)) \Big)\\
&+ \frac{1}{2}\Big(ds (X) \g{Y, Z} - ds (Y) \g{X, Z}\Big).
\end{split} \end{equation}
Since $\rho$ closed \cite{Dr}, we have
\begin{equation}\label{part1}
(\nabla _X \rho)(Y, JZ) - (\nabla _Y \rho)(X, JZ) = -(\nabla_{JZ}\rho)(X,Y).
\end{equation}
As an algebraic object, $\nabla_X J$ is a skew-symmetric endomorphism of
$TM$, associated (by $g$-duality) to the section $\nabla_X \omega$ of the
bundle of $J$-anti-invariant 2-forms; it then anti-commutes with $J$, and
commutes with any skew-symmetric endomorphism associated to a section of
$\Lambda^-M$; in particular, $\nabla_X J$ commutes with the endomorphism
corresponding to $\rho_0$ via the metric (which will be still denoted by
$\rho_0$). We thus obtain
\begin{equation*}\notag\begin{split}
\rho(Y, (\nabla_X J)(Z)) - \rho(X, (\nabla_Y J)(Z) )&=
\frac{3s}{2}\Big((\nabla_Y \omega)(X, JZ) -(\nabla_X \omega)(Y, JZ) \Big)\\
&+(\nabla_X \omega)(Y, \rho_0(Z)) -(\nabla_Y \omega)(X, \rho_0(Z)).
\end{split}\end{equation*}
By using the closedness of $\omega$ we derive
\begin{equation}\label{part2}
\begin{split}
\rho(Y, (\nabla_X J)(Z)) - \rho(X, (\nabla_Y J)(Z))  = &
\frac{3s}{2}(\nabla_{JZ} \omega) (X,Y) \\
&- (\nabla_{{\rho_0} (Z)} \omega)(X,Y).
\end{split}
\end{equation}
Substituting  (\ref{part1}) and (\ref{part2}) in (\ref{basicpart}), we
finally get
\begin{equation}\label{narho0}
\nabla _Z \rho_0  = -\frac{3}{2}ds(Z)\omega - 2 C(JZ)
- \nabla_{\Ric_0(Z)}\omega + ds\wedge JZ^{\flat},
\end{equation}
The $\Lambda^-M$-component  of (\ref{narho0}) gives the identity
(\ref{eq:bianchirho}); the last part of the lemma is immediate.
\end{proof}

It follows that Lemma~\ref{lemma1bis} and hence Proposition \ref{proprho0}
remain true for weakly selfdual K{\"a}hler surfaces, so that on the open set
$M_0$ where $\rho_0 \neq 0$ the almost Hermitian structure $({\bar g}=
\lambda^{-2}g, I)$ defined on $M_0$ by $\rho_0 = \lambda \omega_I$ (see
Lemma \ref{lem:twistor}) is K{\"a}hler.

The theory of \biham/ $2$-forms $\varphi$ does not extend automatically to
the almost K{\"a}hler case: there is no reason, in general, to suppose that
the trace and pfaffian of $\varphi$ are Poisson-commuting holomorphic
potentials, nor can we appeal to the open mapping theorem when $J$ is not
integrable. On the other hand Proposition~\ref{prop:dxi-deta} does
generalize in the following sense: if we write $\sigma=\xi+\eta$ and
$\pi=\xi\eta$, then $d\xi$ and $d\eta$ are orthogonal and
$Jd\sigma=2Id\lambda$ on the closure of $M_0$.

Fortunately, when $\varphi$ is the Ricci form, we can show more.

\begin{lemma}\label{ak-killing}
On a weakly selfdual almost K{\"a}hler $4$-manifold with $J$-invariant
Ricci tensor, $K=J{\rm grad}_g s$ is a Killing vector field.
\end{lemma}
\begin{proof}
On $M_0$, the Ricci tensors of both ${\bar g}= \lambda^{-2}g$ and $g$ are
$I$-invariant, and therefore $I {\rm grad}_g\lambda$ is Killing vector field
(with respect to both metrics)~\cite{AG1}. Since $Ids=2Jd\lambda$,
$Jds=2Id\lambda$ and $K$ is a Killing vector field on $M_0$. On the other
hand if $\lambda$ vanishes identically on an open set $U$ then $g$ is
Einstein on $U$, so that $s$ is constant, and $K$ is a trivial Killing vector
field. Hence by continuity $K$ is a Killing vector field everywhere.
\end{proof}
Because of this observation, it is natural to strengthen the definition of
\biham/ in the almost K{\"a}hler case: we say that a closed $J$-invariant
$2$-form $\varphi$ on an almost K{\"a}hler $4$-manifold is {\it \monoham/} if
its trace-free part $\varphi_0$ is a twistor $2$-form and its trace $\sigma$
is a momentum map for a Killing vector field.

Note that for any Killing vector field $K$, $\nabla_X(\nabla K) = R^g_{K,X}$
and so the $1$-jet $\{K,\nabla K \}$ is parallel with respect to a globally
defined connection, cf. \cite{kostant}.  Hence if $K$ vanishes on an open
set so does $\nabla K$, and therefore $K$ vanishes on any connected
component meeting that open set.  It follows that if $\varphi$ is \monoham/,
the open set $M_0$ where $\varphi_0\neq0$ is dense or empty in each
connected component.

In particular, we obtain the following generalization of
Proposition~\ref{prop:tri}.

\begin{prop}\label{ak-prop:tri} Let $(M,g,J,\omega)$ be a connected
weakly selfdual almost K{\"a}hler $4$-manifold with $J$-invariant Ricci
tensor. Then one of the following holds:
\begin{enumerate}
\item $\rho_0$ is identically zero; then, $(g,J,\omega)$ is
an  Einstein almost K{\"a}hler $4$-manifold; or

\item the scalar curvature $s$ of $g$ is constant, but $\rho_0$ is not
identically zero; then, $(g,J)$ is obtained from a K{\"a}hler surface
$(g,I)$ with two distinct constant principal Ricci curvatures, $\lambda$ and
$\mu$, in the following manner: $J$ equals to $I$ on the
$\lambda$-eigenspace of the Ricci tensor, but to $-I$ on the
$\mu$-eigenspace; hence $I$ is compatible with the opposite orientation of
$(M,J)$; or

\item $s$ is not constant and $g$ is selfdual; or

\item $W^-$ and $\rho_0$ have no zero; then,  $({\bar g} =
\lambda^{-2}g,I)$ is a globally defined extremal K{\"a}hler metric of
non-constant scalar curvature, which is compatible with
the opposite orientation of $(M,J)$.
\end{enumerate}
\end{prop}

\subsection{Weak selfduality, \monoham/ $2$-forms, and a conjecture of
Goldberg}

The existence of {\it non-integrable} almost K{\"a}hler $4$-manifolds listed
in (i)-(iv) of Proposition \ref{ak-prop:tri} appears to be a non-trivial
problem.  We collect below some remarks and known results on this issue:

\begin{bulletlist}
\item A long-standing conjecture of Goldberg \cite{Go,sekigawa} states that
a compact {Einstein} almost K{\"a}hler manifold must be a
K{\"a}hler-Einstein manifold.  The first {\it local} examples of
non-integrable Einstein almost K{\"a}hler $4$-manifolds have been recently
discovered in \cite{NuP,Arm}; we shall provide new examples (see Example 2
below), but for all these examples the Ricci tensor and the anti-selfdual
Weyl tensor identically vanish.

\item The almost K{\"a}hler $4$-manifolds described in Proposition
\ref{ak-prop:tri}(ii) have been recently studied in \cite{ADM}. It is known
that there are essentially two examples of {\it homogeneous} K{\"a}hler
surfaces $(M,g,I)$ which give rise to homogeneous non-integrable almost
K{\"a}hler $4$-manifolds described in (ii) of Proposition \ref{ak-prop:tri};
specifically, $(M,g,I)$ is either isomorphic to $(SU(2)\ltimes Sol_2)/U(1)$
(in the case when the signature of the hermitian form is (2,2),
cf.~\cite{shima}), or to $({\rm Isom}({\mathbb E}^2)\ltimes Sol_2)/SO(2)$
(in the case when the signature of the hermitian form is (0,2),
cf.~\cite{wall}); see also Example 3 below.  However, there are also many
non-homogeneous examples \cite{AAD}.

\item We believe that there should exist non-integrable examples of almost
K{\"a}hler $4$-manifolds described in Proposition \ref{ak-prop:tri}(iii) and
(iv). However, as we discuss below, it is not clear how to generalize the
constructions of Sections 2 and 3 to obtain such examples.
\end{bulletlist}

According to Proposition \ref{prop:biextremal}, for a weakly selfdual
K{\"a}hler surface, not only the scalar curvature, but also the pfaffian of
the normalized Ricci form is a momentum map for a Killing vector field.
Indeed, this much holds for the trace and the pfaffian of the normalized
$2$-form associated to any \biham/ $2$-form, by
Proposition~\ref{prop:biham}.

Hence one approach to generalize the constructions of Sections 2 and 3 to
the almost K{\"a}hler case is to study \monoham/ $2$-forms such that the
pfaffian is also a momentum map for a Killing vector field. The following
Lemma shows that there are no non-integrable examples with linearly
independent Killing vector fields, which means that there is no direct
generalization of the constructions of Section 2.

\begin{lemma}\label{ak-lem3} For an almost K{\"a}hler $4$-manifold
with a \monoham/ $2$-form $\varphi$, the pfaffian $\pi$ of associated
normalized $2$-form $\tilde\varphi$ is a momentum map for a Killing vector
field $\tilde K=J{\rm grad}_g\pi$ if and only if either $(g,J)$ is a
K{\"a}hler or $\tilde K$ is a constant multiple of $K=J{\rm grad}_g\sigma$.
\end{lemma}
\begin{proof} Since $\pi=\frac14\sigma^2-\lambda^2$ we have
$d\pi=\frac12\sigma d\sigma+\varphi_0(Jd\sigma)$. Straightforward
calculation gives
\begin{equation*}
\nabla(Jd\pi)=d\sigma\wedge Jd\sigma-\frac12 |d\sigma|^2\omega
+\frac12\sigma\nabla(Jd\sigma)-\varphi_0\circ\nabla d\sigma.
\end{equation*}
Since $K=J {\rm grad}_g\sigma$ is Killing by assumption, it follows that
$J{\rm grad}_g\pi$ is Killing if and only if $\varphi_0\circ\nabla d\sigma$
is skew. This is automatic on the open set where $\varphi_0$ vanishes, where
$d\sigma=0$ and hence $d\pi=0$. Therefore we can assume $\varphi_0$ is
nonvanishing and write $\varphi_0=\lambda\omega_I$, where $I$ is a complex
structure of the opposite orientation to $J$.

Now $I\circ \nabla d\sigma$ is skew if and only if $\nabla d\sigma$ is
$I$-invariant; since $J$ and $I$ commute, this means that the
$J$-anti-invariant part of $\nabla d\sigma$ must be $I$-invariant. But $K=J
{\rm grad}_g \sigma$ is Killing, so the $J$-anti-invariant part of
$(\nabla_X d\sigma)(Y)$ is equal to
$$2(\nabla _X \omega)(K,Y) + (\nabla_K \omega)(X,Y).$$
The latter is $I$-invariant if and only if for any vector field $X$ we
have
\begin{equation}\label{ak-p-killing}
(\nabla_{IX} J)(K) = I(\nabla_X J)(K).
\end{equation}
Suppose now that $\nabla_X J =A\neq 0$ for a vector $X$ at some point.
Since $\nabla_{IX} J$, like $A$, is a $J$-anti-invariant endomorphism, we
can write $\nabla_{IX} J = b A+ c JA$ for some
$b,c\in\R$. Equation~\eqref{ak-p-killing} now reads
$$b A(K) + c J A(K) = I A(K).$$
Since $A$ commutes with $I$ and anticommutes with $J$, by applying $A$ to
the both sides we  obtain: $-b K + c JK = - IK$. However $K$ is orthogonal
to both $JK$ and $IK$, and so $b=0$ and $c=\pm 1$. Thus, on the open
set where $J$ is non-integrable and $\varphi_0\neq 0$, we have
$J d\sigma=\pm I d\sigma= \pm 2J d\lambda$, so $d\pi$ and $d\sigma$
are linearly dependent.
\end{proof}

The Calabi construction in Section 3 does generalize to the almost K{\"a}hler
case and generates some new examples of selfdual Ricci-flat almost K{\"a}hler
$4$-manifolds. However, we shall see that there are no non-integrable
examples of non-constant scalar curvature.

\begin{prop}\label{ak-lem2}  Let $(M,g,J,\omega)$ be an almost K{\"a}hler
$4$-manifold with $J$-invariant Ricci tensor and a non-vanishing hamiltonian
Killing vector field $K$.  Suppose that the pair $(\bar g=\lambda^{-2}g,I)$ is
K{\"a}hler, where $\lambda$ is a momentum map for a nonzero multiple of $K$,
and $I$ is equal to $J$ on ${\rm span} (K,JK)$, but to $-J$ on the
orthogonal complement of ${\rm span} (K,JK)$.

Then either $J$ is integrable, or $(g,\omega)$ is given explicitly by
\begin{align}\label{ak-metric}
g &= \frac{W}{z}(z^2 g_{\Sigma} + dz^2) + \frac{z}{W}\Bigl(
dt + \frac{V}{z} dz+\beta\Bigr)^2,\\
\omega &= z W\omega_{\Sigma} + dz\wedge\Bigl(
dt + \frac{V}{z} dz+\beta\Bigr),
\label{ak-omega}\end{align}
where $g_\Sigma$ is a metric on $2$-manifold $\Sigma$ with area form
$\omega_\Sigma$, $\beta$ is a $1$-form on $\Sigma$ with
$d\beta=W\,\omega_{\Sigma}$, and $V+iW$ is an arbitrary holomorphic
function on $\Sigma$.

Conversely any such metric satisfies the above hypotheses, and $J$ is
integrable if and only if the function $W$ is constant.
\end{prop}
\begin{proof} Without loss of generality, we take $\lambda=z$ to be
a momentum map for $K$ (with respect to $\omega$). As in
Proposition~\ref{prop:monax}, cf.~LeBrun~\cite{lebrun}, we may introduce
coordinates such that
\begin{align*}
g&= e^u w(dx^2 + dy^2) + w \,dz^2 + w^{-1} (d t  + \alpha)^2\\
\omega&=e^u w \,dx\wedge dy + dz\wedge (d t  + \alpha)\\
\bar\omega&=z^{-2}(-e^u w \,dx\wedge dy + dz\wedge (d t  + \alpha))
\end{align*}
where $\bar\omega$ is the K{\"a}hler form of $I$ with respect to $\bar
g=z^{-2}g$. The integrability of $I$ together with the closedness
of $\omega$ and $\bar\omega$ yields
\begin{align*}
e^uw &= h(x,y)z,\\
d\alpha&=-w_x dy\wedge dz + w_y dx\wedge dz + h(x,y)dx\wedge dy
\end{align*}
from which we obtain the integrability condition
\begin{equation}\label{eq:w-harm}
w_{xx}+w_{yy}=0.
\end{equation}
The Ricci tensor of $g$ is $J$-invariant if and only if
\begin{equation*}
\Bigl(\frac{zu_z-2}{zw}\Bigr)_x =0
\quad\text{and}\quad \Bigl(\frac{zu_z-2}{zw}\Bigr)_y =0
\end{equation*}
so that we can write
\begin{equation*}
zu_z = f(z)\,zw - 2.
\end{equation*}
Since $e^uw=h(x,y)z$, we obtain
$$(zw)_z + \frac{f(z)}{z}(zw)^2 =0.$$
The latter is explicitly integrated, and we get
\begin{equation*}
z\,w(x,y,z) = \frac{1}{F(z) + G(x,y)}
\end{equation*}
for some functions $F(z)$ and $G(x,y)$. By substituting
into~\eqref{eq:w-harm} we discover that either $F$ or $G$ must
be constant.

If $G$ is constant, then $w_x=w_y=0$, i.e., $g$ is of Calabi type and $J$
is integrable.

Consider now the case when $F$ is a constant; then,
$$w= \frac{W}{z} \qquad\text{and}\qquad e^u = z^2 e^U,$$
where $W(x,y)$ is a positive harmonic
function and $U(x,y)$ is an arbitrary function of $(x,y)$; the almost
K{\"a}hler structure $(g,\omega)$ takes the
form~\eqref{ak-metric}--\eqref{ak-omega} where:
\begin{itemize}
\item $g_{\Sigma} = e^U(dx^2 + dy^2)$;
\item $\omega_{\Sigma}=e^U dx\wedge dy$;
\item $W$ is a positive harmonic function on $\Sigma$;
\item $\alpha$ satisfies $d\alpha = -W_y dx\wedge dz/z + W_x dy\wedge dz/z+
W\omega_{\Sigma}$ and we can locally choose $t$ so that $\alpha=V
dz/z+\beta$, where $V$ is a harmonic conjugate of $W$ and
$d\beta=W\sigma_{\Sigma}$.
\end{itemize}
This almost Hermitian structure $(g,J,\omega)$ is almost \ka with
$J$-invariant Ricci tensor, since $w$, $e^u$ and $\alpha$ solve the
required equations.
\end{proof}

One directly calculates the normalized scalar curvature $s$ of the metric
\eqref{ak-metric}: it is given by
$$s = \frac1{6zWe^U}\Bigl(U_{xx} + U_{yy} + 2e^{U}\Bigr).$$

On a weakly selfdual almost K{\"a}hler $4$-manifold with $J$-invariant
Ricci tensor, we have seen that $s$ is a momentum map
for a Killing vector field. However, $s$ cannot be a multiple of $z$ unless
it vanishes. Hence this construction does not yield any
non-integrable weakly selfdual almost \ka metrics of Calabi type with
{\it non-constant} scalar curvature. However, it does provide the
following new examples of
selfdual {\it Ricci-flat} strictly almost \ka $4$-manifolds.

\smallskip\noindent {\bf Example 2.}
Let $(g,J,\omega)$ be given by (\ref{ak-metric}) and suppose that $s=0$;
this means that $U$ is a solution of the Liouville equation, i.e., that
$g_{\Sigma}$ is the standard metric on an open subset of $S^2$,
while $H= W + iV$ is a non-constant holomorphic function on $\Sigma$ with
positive real part.  If we write $z=r$, we see that the metric
$$g = \frac{W}{r} ( dr^2 + r^2 g_{S^2})  + \frac{r}{W} (dt+\alpha)^2$$
is given by applying the Gibbons--Hawking Ansatz using the harmonic function
$W/r$, which is invariant under dilation with weight $-1$ in the sense that
\begin{equation*}
r\frac\partial{\partial r} \Bigl(\frac Wr\Bigr) = -\frac Wr.
\end{equation*}
(In fact this is the natural scaling weight for $W/r$, since it is
the Higgs field of an abelian monopole on $\R^3$.)

This class of Gibbons--Hawking metrics has been studied before
in~\cite{CaPe:sdc} and~\cite{CaTo:emh}. In addition to the triholomorphic
Killing vector field $\frac\partial{\partial t}$, these metrics also admit a
triholomorphic homothetic vector field $r\frac\partial{\partial r}$.
Therefore, by~\cite{GaTo:hms}, the local quotient by
$r\frac\partial{\partial r}$ is a {\it hyperCR Einstein--Weyl space}.  In
this case the quotient Einstein--Weyl structure was obtained explicitly
in~\cite{CaPe:sdc} and is an Einstein-Weyl space {\it with a geodesic
symmetry}.

The reader is referred to these references for more information.  However,
to the best of our knowledge, the observation that these metrics are almost
K{\"a}hler is new. Note that the K{\"a}hler form is not an eigenform of the
Weyl tensor, showing that the solutions are different from the previously
known examples of Nurowski--Przanovski~\cite{NuP} and Tod, which were
obtained by applying the Gibbons--Hawking Ansatz to a translation-invariant
harmonic function.

\smallskip

We next use the rough classification given by Proposition \ref{ak-prop:tri}
to obtain the following partial result which motivates the further study of
compact weakly selfdual almost K{\"a}hler $4$-manifolds with $J$-invariant
Ricci tensor:

\begin{theo}\label{theo6}
Suppose that there exists a compact weakly selfdual almost K{\"a}hler
$4$-manifold $(M,g,J,\omega)$ with $J$-invariant Ricci tensor, for which the
almost complex structure $J$ is not integrable.  Then one of the following
two alternatives holds:
\begin{enumerate}
\item The scalar curvature of $g$ is a negative constant; then $M$ admits an
Einstein, non-integrable almost K{\"a}hler structure; or
\item $(M,g,J,\omega)$ belongs to case \textup{(iv)} of
Proposition~\ref{ak-prop:tri}, and the globally defined K{\"a}hler structure
$(\bar g, I)$ is isomorphic to an extremal K{\"a}hler metric which is not
locally of cohomogeneity one, on a minimal
ruled surface $S = {\mathbb P}({\mathcal O} \oplus L) \to \Sigma_{\bf g}$
with ${\bf g} \ge 1$ and ${\rm deg} \ L >0$.
\end{enumerate}
\end{theo}
\begin{proof}
We inspect the possible compact {\it non-integrable} almost \ka $4$-manifolds
given by (i)--(iv) of Proposition \ref{ak-prop:tri}.

The case of constant scalar curvature is described by Proposition
\ref{ak-prop:tri}(i) \& (ii). Our claim then follows by
\cite{sekigawa} and \cite[Th.2]{ADM}.

Suppose that $s$ is not constant, i.e., that $K = J {\rm grad}_g s$ is a
non-trivial Killing vector field by Lemma~\ref{ak-killing}. Let $x_0\in M$
be a zero of $K$; then, the isotropy subgroup $H(x_0)$ of the connected
group of isometries of $(M,g)$ is a compact group of dimension at least one;
one can therefore take an $S^1$ in $H(x_0)$. Hodge theory implies that on a
compact almost \ka manifold any isometry which is homotopic to the identity
inside the group of diffeomorphisms is a symplectomorphism (see
e.g.~\cite{lichne}); hence the chosen isometric $S^1$-action is {\it
  symplectic} with respect to $\omega$. Since $x_0$ is a fixed point of the
$S^1$-action, we obtain a hamiltonian $S^1$-action on $(M,\omega)$
\cite{McDuff}.  The manifold is then equivariantly (and orientedly)
diffeomorphic to a rational or a ruled complex surface endowed with a
holomorphic circle action \cite{Audin,Ahara-Hattori,Karshon}. Moreover, in
this case $(M,g,J,\omega)$ is given by Proposition \ref{ak-prop:tri} (iii)
or (iv).

Consider first the case (iii). Since $s$ is not constant, the selfdual Weyl
tensor $W^+$ does not vanish \cite[Cor.1]{AD}. By the Chern-Weil formulae,
the signature of $M$ is strictly positive, and therefore $M$ is
diffeomorphic to ${\mathbb CP}^2$ \cite{BPV}.  Combining the results of
\cite{gursky} and \cite{poon}, one sees that on ${\mathbb CP}^2$ the only
selfdual conformal structure with non-trivial (conformal) Killing vector
field is the standard one. Thus, modulo diffeomorphisms, we may assume that
$g$ is conformal to the standard \ka metric $(g_0,\omega_0)$. Since $\omega$
and $\omega_0$ are both harmonic selfdual 2-forms on $({\mathbb CP}^2,
g_0)$, and since $b^+({\mathbb CP}^2)$=1, we conclude that $\omega = const.
\omega_0$, showing that $J$ is integrable, a contradiction.

Suppose $(M,g,J, \omega)$ is as in Proposition \ref{ak-prop:tri}(iv). Now
$\rho_0$ determines an integrable almost complex structure $I$ compatible
with $g$ and with the opposite orientation of $(M,J)$, such that $({\bar
g}=\lambda^{-2}g, I)$ is an {\it extremal} \ka metric with $I{\rm
grad}_{\bar g} {\bar s} = const.K$. Denote by ${\overline M}$ the smooth
manifold $M$ endowed with the orientation induced by $I$. Thus, the oriented
smooth $4$-manifolds $M$ and ${\overline M}$ both admit complex
structures. Since $M$ is the underlying smooth manifold of a rational or a
ruled complex surface, we conclude as in the proof of Theorem {\ref{theo3}}
that the complex surface $({\overline M}, I)$ is a ruled surface of the form
${\mathbb P}(E)$, where $E \to \Sigma_{\bf g} $ is a holomorphic rank 2
bundle over a compact Riemann surface $\Sigma_{\bf g}$ of genus ${\bf g}$,
and which splits as $E ={\mathcal O} \oplus L$ for a holomorphic line bundle
$L$ of degree ${\rm deg} \ (L) > 0$.

We have to prove that ${\bf g}\ge 1$. Indeed, if ${\bf g}=0$, we obtain the
{Hirzebruch surface} $F_k= {\mathbb P}({\mathcal O} \oplus {\mathcal
O}(-k))$, where $k$ is a positive integer.  As we have already observed
in the proof of Theorem \ref{theo3}, the extremal K{\"a}hler metrics on
these surfaces are the Calabi cohomogeneity-one U(2)-metrics, i.e.
${\bar g}$ is given by the Calabi
construction
\eqref{eq:bianchi}--\eqref{eq:Jbianchi}; since $g = const.{\bar s}^2
{\bar g}$, it follows that $g$ is a cohomogeneity-one metric as well and
therefore ${\rm grad}_{\bar g} {\bar s} = const. {\rm grad}_g s$, showing
that $JK = IK$. Since $s$ is not constant, by Lemma \ref{ak-lem3}, $J$ is
integrable on the open dense subset where $ds \neq 0$, hence
everywhere. We thus conclude that ${\bf g} \ge 1$.  Note that the above
local argument applies to any extremal  K{\"a}hler ${\bar g}$ which
is locally of cohomogeneity one, so that the last
part of the theorem also follows.
\end{proof}

We do not have any examples in case (ii) of the above theorem.  Indeed, the
only examples we know of extremal \ka metrics on the minimal ruled surfaces
in (ii) are locally cohomogeneity-one Calabi-type metrics
\cite{christina,christina2}.

\subsection{Almost \ka $4$-manifolds of constant Lagrangian sectional
curvature}

In this section we deduce another global result from Theorem~\ref{theo6}.
An almost K\"ahler $4$-manifold $(M,g,\omega)$ is said to have (pointwise)
\emph{constant Lagrangian sectional curvature} if the sectional curvature of
$g$, at each point of $M$, is constant on the set of Lagrangian $2$-planes
at that point---recall that the latter are the planes $X\wedge Y$ with
$\omega(X,Y)=0$. One can make the same definition for almost K\"ahler
manifolds of any dimension, but for $2n>4$, any almost K\"ahler
$2n$-manifold of constant Lagrangian sectional curvature is in fact K\"ahler
with constant holomorphic sectional curvature \cite{FF,FFK}.  Conversely, it
is easy to see that any complex space form has constant Lagrangian sectional
curvature (see e.g.~\cite{Gan}).

In four dimensions, the situation is more interesting.  A simple local
calculation (cf. e.g.~\cite{AD}) shows that an almost \ka $4$-manifold has
constant Lagrangian sectional curvature if and only if the Ricci tensor is
$J$-invariant, the Weyl tensor is selfdual, and the \ka form is one of its
roots (i.e., $M$ has \emph{Hermitian Weyl tensor} in the sense
of~\cite{AA}).

The following homogeneous example shows that the integrability for almost
\ka $4$-manifolds with constant Lagrangian sectional curvature does not
follow locally (nor even for complete metrics).

\smallskip \noindent {\bf Example 3.} Consider the homogeneous \ka surface
$M=(SU(2)\ltimes Sol_2)/U(1)$, where $Sol_2$ denotes the real
two-dimensional solvable subgroup of upper triangular matrices in
$SL_2({\mathbb R})$ .

We take the unique left-invariant \ka structure $(g,I)$ on $M$, determined
by the property that the constant principal Ricci curvatures are equal to
$(-1,+1)$, cf.  \cite{shima}.  According to Proposition
\ref{ak-prop:tri}(ii), $(M,g)$ admits an almost \ka structure $J$ with
$J$-invariant Ricci tensor. Since the scalar curvature of $g$ is zero, $g$
is selfdual (with the orientation opposite to $I$)~\cite{gauduchon};
furthermore, $J$ is not integrable \cite{ADM}, and by using the general
formulae in \cite[Ch.7]{besse} one easily checks that $(M,g,J)$ has constant
Lagrangian sectional curvature.

\smallskip
In contrast to this example, there were a number of reasons \cite{AD,AA}
to believe that a {\it compact} almost K\"ahler $4$-manifold of constant
Lagrangian curvature must be a selfdual K{\"a}hler metric. The
conjectured integrability of the almost complex structure has been
proved under some additional assumptions on curvature \cite{AD} or the
topology \cite{AA} of the manifold, but the general question was left open.
As a consequence of Theorem \ref{theo6} we are now able to give a positive
answer.

\begin{cor}\label{theo7}
A compact, $4$-dimensional, almost \ka manifold has constant Lagrangian
sectional curvature if and only if it is a K{\"a}hler selfdual surface.
\end{cor}
\begin{proof} Suppose $(M,g,J,\omega)$ is a compact  almost \ka
$4$-manifold of constant Lagrangian sectional curvature, but for which $J$ is
not integrable. According to \cite[Th.2]{AD} the scalar curvature $s$ of $g$
is not constant; then, by Theorem \ref{theo6}, the smooth manifold
${M}$ is diffeomorphic to a minimal ruled surface. Since any such
surface admits an orientation reversing involution, we conclude that $M$
carries a complex structure which is compatible with the
orientation induced by $\omega$. Then, by \cite[Cor.2]{AA}, the almost
complex structure $J$ must be integrable, contradicting our assumption.
\end{proof}


\begin{thebibliography}{99}


\bibitem{Abreu} M. Abreu,
{\it K{\"a}hler geometry of toric varieties and extremal metrics},
Int. J. Math. {\bf 9}  (1998), 641--651.

\bibitem{Ahara-Hattori} K. Ahara and A. Hattori, {\it Four dimensional
symplectic $S^1$-manifolds admitting moment map}, J. Fac. Sci. Univ.
Tokyo, Sect. IA, 38 (1991), 251--298.

\bibitem{AA} V. Apostolov and J. Armstrong,
{\it Symplectic $4$-manifolds with Hermitian Weyl tensor},
{Trans. Amer. Math. Soc.} {\bf 352} (2000), 4501--4513.

\bibitem{AAD} V. Apostolov, J. Armstrong, and T. Draghici,
{\it Local rigidity of certain classes of almost \ka $4$-manifolds},
preprint 1999, available at arXiv: math.DG/9911197.

\bibitem{AD} V. Apostolov and T. Draghici,
{\it Almost K{\"a}hler $4$-manifolds with $J$-invariant Ricci tensor and
special Weyl tensor}, Q. J. Math. {\bf 51} (2000),
275--294.

\bibitem{ADM} V. Apostolov, T. Draghici, and A. Moroianu,
{\it A splitting theorem for K{\"a}hler manifolds with constant
eigenvalues of the Ricci tensor}, to appear in Int. J. Math., available
at arXiv: math.DG/0007122.

\bibitem{AG1} V. Apostolov and P. Gauduchon,
{\it The Riemannian Goldberg-Sachs Theorem}, Int. J. Math. {\bf 8}
(1997), 421--439.

\bibitem{AG2} V. Apostolov and P. Gauduchon,
{\it Self-dual Einstein Hermitian $4$-manifolds}, preprint 2000,
available at arXiv:mathDG/0003162.

\bibitem{Arm} J. Armstrong,
{\it An Ansatz for Almost-K{\"a}hler, Einstein $4$-manifolds}, preprint
1997, to appear in J. reine angew. Math.

\bibitem{Audin} M. Audin, {\it Hamiltoniens p{\'e}riodiques sur les
vari{\'e}t{\'e}s symplectiques de dimension $4$}, in ``G{\'e}om{\'e}trie
symplectique et m{\'e}canique'', Proceedings 1988 (ed. C. Albert), Springer
Lecture Notes in Math. 1416 (1990).


\bibitem{BPV} W. Barth, C. Peters and A. Van de Ven,
{\it Compact complex surfaces}, Springer-Verlag,
Berlin Heidelberg New York Tokyo, 1984.

\bibitem{besse} A. L. Besse, {\it Einstein manifolds},
Ergeb. Math. Grenzgeb. {\bf 3}, Springer-Verlag, Berlin, Heidelberg,
New York, 1987.

\bibitem{bryant} R. Bryant, {\it Bochner-K{\"a}hler metrics},
J. Amer. Math. Soc. {\bf 14} (2001) 623--715.

\bibitem{calabi1} E. Calabi, {\it Extremal K{\"a}hler metrics},
Seminar on Differential Geometry, Princeton Univ. Press (1982), 259--290.

\bibitem{calabi2} E. Calabi, {\it Extremal K{\"a}hler metrics II},
in: Differential Geometry and Complex Analysis (eds. I. Chavel and
H.M. Farkas), Springer-Verlag, 1985.

\bibitem{CaPe:sdc} D.~M.~J. Calderbank and H.~Pedersen,
{\it Selfdual spaces with complex structures, {E}instein-{W}eyl geometry
and geodesics}, Ann. Inst. Fourier {\bf 50} (2000) 921--963.

\bibitem{CaPe:sdt2} D.~M.~J. Calderbank and H.~Pedersen,
{\it Selfdual Einstein metrics with torus symmetry}, Preprint,
University of Edinburgh.

\bibitem{CaTo:emh} D.~M.~J. Calderbank and K.~P. Tod,
{\it {E}instein metrics, hypercomplex structures and the {T}oda field
equation}, to appear in Diff. Geom Appl. (2001).

\bibitem{kosniowski} J. Carrel, A. Howard and C. Kosniowski,
{\it Holomorphic vector fields on compact complex manifolds},
Math. Ann. {\bf 204} (1973), 73--81.

\bibitem{BYC} B.-Y. Chen, {\it Some topological obstructions to
Bochner-K{\"a}hler metrics and their applications},
J. Differential Geom. {\bf 13} (1978), 574--588.

\bibitem{derdzinski} A. Derdzi\'nski,
{\it Self-dual K{\"a}hler manifolds and Einstein manifolds of dimension four},
Compositio Math. {\bf 49}  (1983), 405--433.

\bibitem{Dr} T. Draghici,
{\it On some $4$-dimensional almost K{\"a}hler manifolds}, Kodai
Math. J. {\bf 18} (1995), 156--163.

\bibitem{FF} M. Falcitelli and  A. Farinola,
{\it Six-dimensional almost \ka manifolds with pointwise constant
antiholomorphic sectional curvature}, Riv. Mat. Univ. Parma {\bf (5)
6} (1997), 143--156.

\bibitem{FFK} M. Falcitelli, A. Farinola, and  O. Kassabov,
{\it Almost \ka manifolds whose antiholomorphic sectional curvature is
pointwise constant}, Rend. Mat. Appl. Ser. {\bf VII}, {\bf 18} (1998),
151--166.

\bibitem{Gan} G. Ganchev,
{\it On Bochner curvature tensors in almost Hermitian manifolds},
Pliska Studia Matematica Bulgarica {\bf 9} (1987), 33--43.

\bibitem{gauduchon} P. Gauduchon, {\it Surfaces k{\"a}hl{\'e}riennes dont la
courbure v{\'e}rifie certaines conditions de positivit{\'e}}, in
``G{\'e}om{\'e}trie riemannienne en dimension 4'', S{\'e}m. A. Besse, 1978/79.

\bibitem{GaTo:hms} P.~Gauduchon and K.~P. Tod,
{\it Hyperhermitian metrics with symmetry},
J.~Geom. Phys. {\bf 25} (1998) 291--304.

\bibitem{Go} S. I. Goldberg,
{\it Integrability of almost K{\"a}hler manifolds},
Proc. Amer. Math. Soc. {\bf 21} (1969), 96--100.

\bibitem{G:kstv} V.~Guillemin,
{\it K{\"a}hler structures on toric varieties},
J. Differential Geom.  {\bf 40}  (1994),  285--309.

\bibitem{gursky} M. Gursky, {\it Conformal vector fields on
four-manifolds with negative scalar curvature}, Math. Z. {\bf 232} (1999),
265--273.

\bibitem{hwang-maschler} A. Hwang and G. Maschler,
{\it Central K{\"a}hler metrics with non-constant central curvature},
Preprint 2000.

\bibitem{jelonek} W. Jelonek,
{\it Compact K{\"a}hler surfaces with $\delta W^-=0$}, Preprint 2000.

\bibitem{jelonek2} W. Jelonek,
{\it Extremal K{\"a}hler AC-surfaces}, Preprint 2001.

\bibitem{Karshon} Y. Karshon, {\it Periodic Hamiltonian flows on
four-dimensional manifolds} in ``Contact and symplectic geometry'',
(Cambridge, 1994), 43--47, Publ. Newton Inst., 8, Cambridge Univ. Press,
Cambridge, 1996.

\bibitem{kostant} B. Kostant, {\it Holonomy and the Lie algebra of
    infinitesimal motions of a Riemann manifold},
    Trans. Amer. Math. Soc. {\bf 80} (1955), 528--542.

\bibitem{kotschick} D. Kotschick,
{\it Orientations and geometrisations of compact complex surfaces},
Bull. London Math. Soc. {\bf 29}, No. {\bf 2} (1997), 145--149.

\bibitem{lebrun} C. R. LeBrun, {\it Explicit self-dual metrics on ${\C P}^2
\# \cdots\# {\C P}^2$}, J. Differential Geom. {\bf 34} (1991), 223--253.


\bibitem{lebrun-singer} C. R. LeBrun and M. Singer,
{\it Existence and deformation theory of scalar flat K{\"a}hler
metrics on compact complex surfaces}, Inv. Math. {\bf 112} (1993],
273--313.

\bibitem{lichne} A. Lichnerowicz, {Th{\'e}orie des groupes de
transformations}, Dunod, Paris, 1958.

\bibitem{McDuff} D. McDuff,
{\it The moment map for circle actions on symplectic manifolds,}
J. Geom. Phys. {\bf 5} (1988), 149--160.

\bibitem{maruyama} M. Maruyama,
{\it On automorphism groups of ruled surfaces}, J.  Math. Koyoto
Univ. {\bf 11} (1971), 89--112.

\bibitem{maschler} G. Maschler, {\it Central K{\"a}hler metrics},
Preprint 2000.

\bibitem{mats-tanno} M. Matsumoto and S. Tanno,
{\it On K{\"a}hler spaces with parallel or vanishing Bochner curvature
tensor}, Tensor N. S. {\bf 27} (1973), 291--294.

\bibitem{NuP} P. Nurowski and M. Przanowski,
{\it A four-dimensional example of Ricci flat metric admitting almost
K{\"a}hler non-K{\"a}hler structure}, Classical Quantum Gravity
{\bf 16} (1999), no. {\bf 3}, L9--L13.

\bibitem{pontecorvo} M. Pontecorvo,
{\it On twistor spaces of anti-self-dual Hermitian surfaces},
Trans. Amer. Math. Soc. {\bf 331} (1992), 653--661.

\bibitem{poon} Y.-S. Poon,
{\it Compact self-dual manifolds with positive scalar curvature},
J. Differential Geom. {\bf 24} (1986), 97--132.

\bibitem{sekigawa} K. Sekigawa,
{\it On some compact Einstein almost K{\"a}hler manifolds},
J. Math. Soc. Japan {\bf 36} (1987), 677--684.

\bibitem{shima} H. Shima,
{\it On homogeneous K{\"a}hler manifolds with non-degenerate canonical
Hermitian form of signature $(2, 2(n-1))$}, Osaka J. Math. {\bf 10}
(1973), 477--493.

\bibitem{christina} C. T{\o}nnesen-Friedman,
{\it Extremal K{\"a}hler metrics on minimal ruled surfaces}, J. reine
angew. Math. {\bf 502} (1998), 175--197.

\bibitem{christina2} C. T{\o}nnesen-Friedman,
{\it Extremal Kahler metrics and Hamiltonian functions II},
to appear in Glasgow Math. Journal.

\bibitem{wall} C. T. C. Wall,
{\it Geometric structures on compact complex analytic surfaces},
Topology {\bf 25} (1986), 119--153.

\end{thebibliography}
\end{document}